\documentstyle[11pt]{article}

\newtheorem{eg}{{\bf Example}}[section]

\newtheorem{theo}{{\bf Theorem}}[section]
\newtheorem{prop}{{\bf Proposition}}[section]
\newtheorem{cor}[theo]{{\bf Corollary}}

\font\bbb=msbm10 scaled\magstep1

\newcommand{\CC}{\mbox{\bbb C}}
\newcommand{\FF}{\mbox{\bbb F}}
\newcommand{\NN}{\mbox{\bbb N}}

\newcommand{\RR}{\mbox{\bbb R}}
\newcommand{\ZZ}{\mbox{\bbb Z}}

\def\full{~\, {\rm \kern.22em
 \vrule width.04em
    height1.2ex depth-.05ex
 \kern-.6em \in}\,}

\newcommand{\AAA}{A^{\!\!\!^{^{\circ}}}}

\newcommand{\CCC}{C^{\!\!\!^{^{\circ}}}}
\newcommand{\gammaint}{\gamma^{\!\!\!^{{\circ}}}}
\newcommand{\sigmaint}{\sigma^{\!\!\!^{{\circ}}}}

\newcommand{\TTT}{T^{\!\!\!^{^{\circ}}}}

\newcommand{\vvv}{v^{\!\!\!^{_{\circ}}}}
\newcommand{\bdA}{A^{\!\!\!^{^{\bullet}}}}
\newcommand{\bdB}{B^{\!\!\!^{^{\bullet}}}}
\newcommand{\bdC}{C^{\!\!\!^{^{\bullet}}}}

\newcommand{\bdv}{v^{\!\!\!^{_{\bullet}}}}
\newcommand{\bdsi}{\sigma^{\!\!\!^{_{\bullet}}}}
\newcommand{\bdDelta}{{\Delta^{\!\!\!\!^{^{\bullet}}}}{}^{d+1}}
\newcommand{\TPSS}{S^{\hspace{.2mm}2}\! \times \hspace{-3.4mm}_{-} \, S^1}
\newcommand{\TPSSD}{S^{\hspace{.2mm}d-1}\! \times \hspace{-3.8mm}_{-} \, S^1}

\begin{document}



\title{\bf Minimal Triangulations of Manifolds}
\author{{\Large \bf Basudeb Datta} \\[2mm]
Department of Mathematics \\
Indian Institute of Science \\
Bangalore 560\,012, India. \\[1mm]
{\it dattab@math.iisc.ernet.in} }
\date{January 2007}

\maketitle






\vspace{-4mm}

\hrule


\section*{Abstract}


In this survey article, we are interested on minimal
triangulations of closed pl manifolds. We present a brief survey
on the works done in last 25 years on the following: (i) Finding
the minimal number of vertices required to triangulate a given
pl  manifold. (ii) Given positive integers $n$ and $d$,
construction of $n$-vertex triangulations of different
$d$-dimensional pl manifolds. (iii) Classifications of all  the
triangulations of a given pl manifold with same number of
vertices.

In Section 1, we have given all the definitions which are required
for the remaining part of this article. In Section 2, we have
presented a very brief history of triangulations of manifolds. In
Section 3, we have presented examples of several vertex-minimal
triangulations. In Section 4, we have presented some interesting
results on triangulations of manifolds. In particular, we have
stated the Lower Bound Theorem and the Upper Bound Theorem. In
Section 5, we have stated several results on minimal
triangulations without proofs. Proofs are available in the
references mentioned there.

\bigskip


\hrule




\section{Preliminaries}

\begin{description}
\item[Affine Subspaces of $\RR^n$ and Linear Maps.]
The space $\{(x_1, \dots, x_n) : x_i \mbox{ is real for } 1\leq
i\leq n\}$ is denoted by $\RR^n$. For us, $\RR^n_+= \{(x_1,
\dots, x_n)\in \RR^n : x_n\geq 0\}$, $I= [-1, 1]\subseteq \RR$ and
$\NN=\{0, 1, 2, \dots\}\subseteq \{0, \pm1, \pm2, \dots\}=\ZZ$.

An {\em affine} subspace $V \subseteq\RR^m$ (of dimension $n$) is
a translated vector subspace (of dimension $n$). So, $V \subseteq
\RR^m$ is an affine subspace if $a_1, \dots, a_r\in V$ and
$\lambda_1, \dots, \lambda_r \in [0, 1]$ with
$\sum_{i=1}^r\lambda_i=1$ implies $\sum_{i=1}^r\lambda_i a_i\in
V$.

A map $f \colon V\to \RR^d$, from an affine subspace $V$ of
$\RR^m$, is called (affine) {\em linear} if $f(\sum_{i = 1}^r
\lambda_i a_i) = \sum_{i = 1}^r \lambda_i f(a_i)$.

Clearly, if $V \subseteq \RR^m$ is an $(m-1)$-dimensional affine
subspace then $\RR^m\setminus V$ has two connected components,
say $H_1$ and $H_2$. The subsets $V\cup H_1$ and $V\cup H_2$ are
called the (closed) {\em half-spaces} determined by $V$.

If the smallest affine subspace in $\RR^m$ containing $n$ points
$v_1, \dots, v_n$ is $(n - 1)$ dimensional (equivalently, $v_2 -
v_1, \dots, v_n - v_1$ are linearly independent), then we say
that the points $v_1, \dots, v_n$ in $\RR^{m}$ are {\em affinely
independent}.

\item[Joins and Cones.]
If $A$, $B$ are subsets of $\RR^n$, then their {\em join} $AB$ is
the subset $\{\lambda a + \mu b : a \in A, b \in B, \lambda, \mu
\in [0, 1]$ and $\lambda + \mu = 1\}$. So, $AB$ consists of all
points on line segments (arcs) with endpoints in each of $A$ and
$B$. If $A=\emptyset$ then we define $AB=B$. If $A=\{a\}$ then
$AB$ is also denoted by $aB$. A join $aB$ is called a {\em cone}
(with vertex $a$ and base $B$) if $a\notin B$ and $b_1$, $b_2\in
B$, $b_1\not=b_2$ then $ab_1\cap ab_2=\{a\}$.

\item[Polytopes and Simplices.]
A subset $C\subseteq \RR^m$ is called {\em convex} if for each
pair of points $a, b\in C$ the arc $ab\subseteq C$. For a set $A$
(possibly empty) in $\RR^m$, the smallest convex set containing
$A$ is called the {\em convex hull of $A$} and is denoted by
$\langle A\rangle$. A {\em polytope} is a convex hull of a finite
set. A polytope $C$ is said to be an $n$ {\em dimensional}
polytope (or {\em $n$-polytope\,}) if the smallest affine
subspace containing $C$ is $n$ dimensional. By convention, the
empty set is a polytope of dimension $-1$.

A point $v$ in a polytope $C$ is called a {\em vertex} if $v \in$
arc $ab \subset C$ implies $v$ is $a$ or $b$.

Clearly, an $n$-polytope has at least $n+1$ vertices. If an
$n$-polytope has exactly $n + 1$ vertices then it is also called
an {\em $n$-simplex}. So, $\langle\{v_0, v_1, \dots,
v_n\}\rangle$ is an $n$-simplex if and only if $v_0, v_1, \dots,
v_n$ are affinely independent. An $n$-simplex with vertices $v_0,
v_1, \dots, v_n$ is denoted by $v_0v_1 \cdots v_n$.

If $A = v_0\cdots v_k$ is a $k$-simplex then $\hat{A} :=
\sum_{i=0}^k\frac{1}{k+1}v_i$ is called the {\em barycentre} of
$A$.

\item[Faces of a Polytope.]
Let $C$ be an $n$-polytope in $\RR^m$. If $V$ is an
$(m-1)$-dimensional affine subspace such that $C$ is in one of the
half-space determined by $V$ then $C\cap V$ is called a {\em
face} of $C$ and is denoted by $C\cap V < C$. Clearly, a face of a
polytope is a polytope. If $\emptyset \neq D< C$ and $D\neq C$
then $D$ is called a {\em proper face} of $C$. The union of all
the proper faces of an $n$-polytope $C$ ($n\geq 1$) is called the
{\em frontier} of $C$ and is denoted by $\bdC$. The subset $\CCC
:= C \setminus \bdC$ is called the {\em interior} of $C$. For a
0-polytope (i.e., for a vertex) $v$ we define $\vvv =v$ and
$\bdv=\emptyset$. So, for a polytope $C$, $C=\CCC\sqcup\bdC$.

\item[Simplicial and Stacked Polytopes.]
A polytope is called {\em simplicial} if its proper faces are
simplices.

A simplicial $d$-polytope $P$ is called {\em stacked} if there is
a sequence $P_1, \dots, P_k$ of simplicial $d$-polytopes such
that $P_1$ is a simplex, $P_k = P$ and $P_{j+1}$ can be
constructed from $P_j$ by attaching a $d$-simplex along a
$(d-1)$-face of $P_j$ for $1\leq j\leq k-1$.

\item[Polyhedra and Subpolyhedra.]
A subset $P\subseteq \RR^n$ is called a {\em polyhedron} if each point
$a$ in $P$ has a cone neighbourhood (in $P$) $N=aL$, where $L$ is compact; $N$ and
$L$ are called the {\em star} and the {\em link} of $a$ in $P$
respectively. We write $N=N_a(P)$ and $L=L_a(P)$. If $P$ and $Q$ are
polyhedra and $Q\subseteq P$ then $Q$ is called a {\em subpolyhedron} of
$P$.
\item[Piecewise-Linear (PL) Maps.] A map $f\colon P\to Q$, where $P$ and $Q$
are polyhedra, is called {\em piecewise-linear} (in short {\em
pl\,}) if each $a \in P$ has a star $N = aL$ such that $f(\lambda
a + \mu x)= \lambda f(a) + \mu f(x)$, for all $x \in L$ and
$\lambda, \mu \geq 0$, $\lambda + \mu = 1$. Moreover, if $f$ is a
homeomorphism then $f$ is called a {\em pl homeomorphism}. A pl
map $f \colon P\to Q$ is called a {\em pl embedding} if $f$ is
injective and $f(P)$ is a subpolyhedron of $Q$. \newline [Check
that $f\colon P \to Q$ a pl homeomorphism implies $f^{-1} \colon
Q \to P$ is pl.]

\item[PL Manifolds.] A polyhedron $M$ is called an $n$-dimensional {\em
pl manifold} (or a {\em pl $n$-manifold}) if each $x \in M$ has a
neighbourhood in $M$ which is pl homeomorphic to an open set in
$\RR^n_+$. The set $\partial M$ consisting of points
corresponding to $\RR^{n-1}\times \{0\}\subseteq\RR^n_+$ is
called the {\em boundary} of $M$. If $\partial M=\emptyset$ then
$M$ is called a {\em pl manifold without boundary}. A compact pl
manifold without boundary is also called a {\em closed pl
manifold}.

[Well-defineness of $\partial M$ follows from the following: {\em
Let $U$ and $V$ be open in $\RR^n_+$ and $f \colon U \to V$ be a
pl homeomorphism. If $x \in \RR^{n - 1} \times \{0\} \cap U$ then
$f(x) \in \RR^{n - 1} \times \{0\}$.}]

Clearly, if $M$ and $N$ are pl manifolds of dimensions $m$ and $n$
respectively then $M \times N$ is a pl $(m+n)$-manifold and
$\partial(M\times N) = (M\times \partial N)\cup (\partial M
\times N)$.

Let $T= ([-2, 2]\times [-2, 2]\setminus (-1, 1)\times (-1,
1))\times [-1, 1]\subseteq\RR^3$. Clearly,  $\TTT$ (the interior
of $T$) $= ((-2, 2)\times (-2, 2)\setminus [-1, 1]\times [-1,
1])\times (-1,1)$. Then $T$ and $T\setminus\TTT$ are pl manifolds
and $\partial T=T\setminus\TTT$. Observe that $\partial T$ is
homeomorphic to the torus $S^{\,1}\times S^{\,1}$.

\item[PL Balls and PL Spheres.] A polyhedron $M$ is called a {\em
pl $n$-ball} if it is pl homeomorphic to $I^{\,n}$. A polyhedron
$M$ is called a {\em pl $n$-sphere} if it is pl homeomorphic to
$\partial I^{\,n+1}$. So, a pl $(n+1)$-ball is a pl $(n +
1)$-manifold having a pl $n$-sphere as boundary. If $C$ is an
$n$-polytope then $C$ is a pl $n$-ball with boundary $\bdC$.

\item[Simplicial Complex.]
A finite collection $K$ of simplices in some $\RR^n$ is called a
{\em simplicial complex} if (i) $\alpha \in K$, $\beta<\alpha$
imply $\beta \in K$ and (ii) $\sigma$, $\gamma\in K$ imply
$\sigma \cap \gamma < \sigma, \gamma$.

For $i=0, 1$, the $i$-simplices in a simplicial complex $K$ are
also called the {\em vertices} and {\em edges} of $K$,
respectively. The set of vertices is called the {\em vertex set}
of $K$ and is denoted by $V(K)$. For a simplicial complex $K$,
the maximum of $k$ such that $K$ has a $k$-simplex is called the
{\em dimension} of $K$.

A simplex $\sigma$ in a simplicial complex $K$ is called {\em
maximal} if $\sigma < \gamma \in K$ implies $\gamma = \sigma$.
Clearly, a simplicial complex is uniquely determined by its
maximal simplices.

A simplicial complex is called {\em pure} if all the maximal
simplices are of same dimension. A maximal simplex in a pure
simplicial complex is also called a {\em facet}.

A simplicial complex of dimension $\leq 1$ is called a {\em
graph}.

\item[Geometric Carrier.] If $K$ is a simplicial complex then
$|K| := \bigcup_{\sigma \in K} \sigma$ is a compact polyhedron
and is called the {\em geometric carrier} of $K$ or the {\em
underlying polyhedron} corresponding to $K$.

\item[Subcomplex.] If $K$ and $L$ are simplicial complexes and
$L \subseteq K$ then $L$ is called a {\em subcomplex} of $K$. We
consider $\emptyset$ to be a subcomplex of every simplicial
complex.

For a simplicial complex $K$, if $U \subseteq V(K)$ then $K[U]$
denotes {\em induced subcomplex} of $K$ on the vertex-set $U$
(i.e., $K[U] = \{\sigma \in K : $ vertices of $\sigma$ are in
$U\}$).

\item[Simplicial Maps.] Let $K$ and $L$ be two simplicial complexes. A map
$f \colon |K| \to |L|$ is called {\em simplicial} if
$f|_{\sigma}$ is linear and $f(\sigma)$ is a simplex of $L$ for
each $\sigma \in K$.

\item[Abstract Simplicial Maps.]
Let $K$ and $L$ be two simplicial complexes. A map $\varphi
\colon V(K) \to V(L)$ is called an {\em abstract simplicial map}
if  $\langle A \rangle$ is a simplex in $K$ implies $\langle
\varphi(A)\rangle$ is a simplex in $L$, for every $A \subseteq
V(K)$.

Let $\varphi \colon K \to L$ be an abstract simplicial map. If
$x\in |K|$ then there exists a unique simplex $\sigma\in K$ such
that $x\in \sigmaint$. Let $\sigma = v_0v_1\cdots v_k$ and $x =
t_0v_0 + \cdots + t_kv_k$ where $t_i \in [0, 1]$ for $0 \leq
i\leq k$ and $t_0+ \cdots + t_k =1$. Define $|\varphi|(x) =
t_0\varphi(v_0) + \cdots + t_k\varphi(v_k)$. This defines a
simplicial map $|\varphi| \colon |K| \to |L|$.

\item[Isomorphisms.]
A bijection $\varphi \colon V(K) \to V(L)$ is called an {\em
isomorphism} if both $\varphi$ and $\varphi^{-1}$ are abstract
simplicial maps. Two simplicial complexes $K$ and $L$ are called
{\em isomorphic} (denoted by $K \cong L$) if such an isomorphism
exists. We identify two simplicial complexes if they are
isomorphic. Clearly, if $\varphi$ is an isomorphism then
$|\varphi|$ is a pl homeomorphism.

An isomorphism from a simplicial complex $K$ to itself is called
an {\em automorphism} of $K$. All the automorphisms of $K$ form a
group under composition, which is denoted by ${\rm Aut}(K)$.

If $\varphi \colon V(K) \to V(L)$ is an isomorphism then define
$\Phi \colon K \to L$ as $\Phi(v_0v_1\cdots v_k) =
\langle\varphi(\{v_0, v_1, \dots, v_k\})\rangle$. Clearly, $\Phi$
is a bijection and $\alpha < \beta$ if and only if $\Phi(\alpha)
< \Phi(\beta)$. Conversely, any such bijection $\Phi \colon K \to
L$ defines an isomorphism $\Phi|_{V(K)}$.

\item[$f$-vector and Euler characteristic.]
If $f_i(K)$ denote the number of $i$-simplices ($0\leq i\leq d$)
in a $d$-dimensional simplicial complex $K$ then $(f_0(K),
f_1(K), \dots$, $f_{d}(K))$ is called the {\em $f$-vector} of $K$
and the number $\chi(K) := \sum_{i=0}^{d} (-1)^i f_i(K)$ is
called the {\em Euler characteristic} of $K$. (Formally we take
$f_{-1} := 1$.)

A simplicial complex $K$ is called {\em $k$-neighbourly} if the
convex hull of any set of $k$ vertices is a $(k-1)$-simplex of $K$
(i.e., $f_{k-1}(X) ={f_0(K) \choose k}$.

\item[Face polynomial and $h$-vector.]
The {\em face polynomial} of a $d$-dimensional simplicial complex
$K$ is
$$
f_K(x) := \sum_{i=-1}^{d} f_i(K)x^{d-i}.
$$

The polynomial $h_K(x) := f_K(x-1)$ is called the {\em
$h$-polynomial} of $K$. The {\em $h$-vector} of $K$ is $(h_0(K),
\dots, h_{d+1}(K))$, where $h_K(x) = \sum_{j=0}^{d+1} h_jx^{d+1-
j}$. Equivalently,
$$
h_j(K) = \sum_{i=-1}^{j-1}(-1)^{j-i-1}{d-i \choose j-i-1}f_i(K)
$$
for $0\leq j\leq d+1$. Observe that $h_{d+1}(K) = (- 1)^{d+1}(1 -
\chi(K))$ and, for $0\leq i\leq d$,
$$
f_{i-1}(K) = \sum_{j=0}^{i}{d+1-j \choose i-j}h_j(K).
$$

\item[Join of Complexes.] Two simplicial complexes $K$ and $L$
(in $\RR^N$) are called {\em independent} if $\alpha\beta$ is an
$(m+n+1)$-simplex for each $m$-simplex $\alpha$ in $K$ and each
$n$-simplex $\beta$ in $L$ for $m, n\geq 0$. If $K$ and $L$ are
independent then we define $K\ast L= K \cup L \cup \{\alpha\beta :
\alpha\in K, \beta\in L\}$. The simplicial complex $K\ast L$ is
called the (simplicial) {\em join} of $K$ and $L$.

If $K$ and $L$ are two simplicial complexes in $\RR^n$ and $\RR^m$
respectively, then we can define their join in a bigger space.
More explicitly,  let $i_1 \colon \RR^n\to \RR^{n + m + 1}$, $i_2
\colon \RR^m \to \RR^{n + m + 1}$ be the maps given by $i_1(x_1,
\dots, x_n)= (x_1, \dots, x_n, 0, \dots, 0)$ and $i_2(x_1, \dots,
x_m)= (0, \dots, 0, x_1, \dots, x_m, 1)$. Let $K_1 : =
\{i_1(\alpha) : \alpha \in K\}$ and $L_1 : = \{i_2(\beta) : \beta
\in L\}$. Then $K_1 \cong K$, $L_1 \cong L$ and $K_1$ and $L_1$
are independent simplicial complexes in $\RR^{n + m + 1}$. We
define $K \ast L = K_1 \ast L_1$.

\item[Stars, Links and Degrees.] Let $K$ be a simplicial complex
and $\gamma \in K$. Let ${\rm st}_K(\gamma)$ be the subcomplex of
$K$ whose maximal simplices are those maximal simplices of $K$
which contain $\gamma$ as a face. This subcomplex is called the
{\em star} of $\gamma$ in $K$.

Let $K$ be a simplicial complex and $\gamma \in K$. Let ${\rm
lk}_K(\gamma) := \{\beta : \beta\cap \gamma = \emptyset,
\beta\gamma \in K\}$. Then, ${\rm lk}_K(\gamma)$ is a subcomplex
of $K$ and is called the {\em link} of $\gamma$ in $K$. The
number of vertices in the link of $\gamma$ in $K$ is called the
{\em degree} of $\gamma$ and is denoted by $\deg_K(\gamma)$.

If $\sigma$ is a simplex in $\RR^n$ then ${\rm Cl}(\sigma) :=
\{\beta : \beta <\sigma\}$ and ${\rm Bd}(\sigma) := \{\beta :
\beta <\sigma, \beta \neq \sigma\}$ are simplicial complexes.
Clearly, $|{\rm Cl}(\sigma)|= \sigma$ and $|{\rm Bd}(\sigma)| =
\bdsi$.

If $\alpha$ is a simplex in a simplicial complex $K$ then ${\rm
Cl}(\sigma)$ and ${\rm lk}_K(\alpha)$ are independent and ${\rm
st}_K(\alpha)$ is the join of ${\rm Cl}(\alpha)$ and ${\rm
lk}_K(\alpha)$.

\item[Subdivisions and Combinatorially Equivalent Complexes.]
A simplicial complex $L$ is called a {\em subdivision} of a
simplicial complex $K$ (denoted by $L\lhd K$) if each simplex in
$L$ is contained in a simplex in $K$ and $|L|= |K|$. Two
simplicial complexes $K$ and $L$ are called {\em combinatorially
equivalent} (denoted by $K \approx L$) if there exist
subdivisions $K^{\prime}\lhd K$ and $L^{\prime}\lhd L$ such that
$K^{\prime} \cong L^{\prime}$. So (by Proposition \ref{prop5}),
$K \approx L$ if and only if $|K|$ and $|L|$ are pl homeomorphic.
Clearly, `$\approx$' is an equivalence relations.

For $\gamma\in K$ and $a\in \gammaint$, consider the simplicial
complex (on the vertex-set $V(K) \cup \{a\}$) $K^{\prime} =
\{\delta \in K : \gamma \not< \delta\} \cup \{a, a\alpha :  \alpha
< \delta$ where $\gamma < \delta \in K$ and $\sigma \neq \alpha
\neq \gamma, \}$. Then $K^{\prime}$ is a subdivision of $K$ and is
called the subdivision obtained from $K$ by {\em starring} at $a$
(or {\em starring the vertex $a$ in $\gamma$}) . We also say that
$K$ is obtained from $K^{\prime}$ by {\em collapsing} the vertex
$a$.

\item[Stellar Subdivisions.] A simplicial complex $K_1$ is called
a {\em stellar} subdivision of $K$ if $K_1$ is obtained from $K$
by starring (successively) at finitely many points. Two complexes
$K$ and $L$ are called {\em stellar equivalent} if they have
isomorphic stellar subdivisions.

Let $A_1, \dots, A_n$ be all the simplices of a simplicial
complex $K$ of dimension $\geq 1$ such that $\dim(A_1) \geq
\dim(A_2) \geq \cdots \geq \dim(A_n)$. Choose $a_i \in \AAA{}_i$
for $1 \leq i \leq n$. Let $K^{(1)}$ be the stellar subdivision of
$K$ obtained by starring at $a_1, \dots, a_n$ successively. Then
$K^{(1)}$ is called a {\em first derived subdivision} of $K$. For
$r\geq 2$, a {\em $r$-th derived subdivision} of $K$ is defined
inductively by $K^{(r)} = (K^{(r-1)})^{(1)}$. If $a_i = \hat{A}_i$
(the barycentre of $A_i$) for each $A_i \in K$ then the first
derived subdivision is called the {\em first barycentric
subdivision} of $K$. Similarly, we can define the {\em $r$-th
barycentric subdivision}.

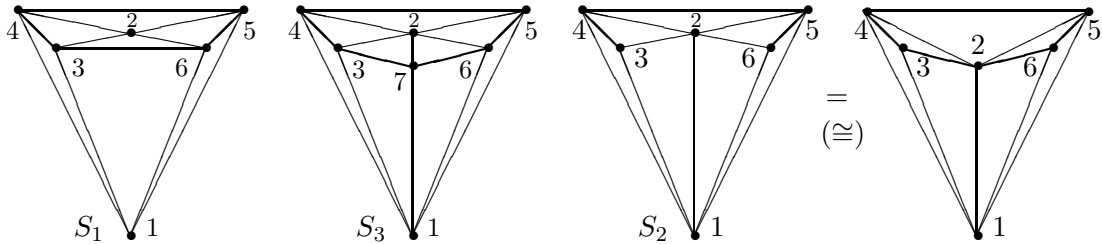
\setlength{\unitlength}{2.5mm}
\begin{picture}(57.2,15)(2.5,-1)



\put(2.77,10){$_{\bullet}$} \put(0.77,12){$_{\bullet}$}
\put(12.77,12){$_{\bullet}$} \put(10.77,10){$_{\bullet}$}
\put(6.77,0){$_{\bullet}$} \put(6.77,10.8){$_{\bullet}$}

\put(6.8,11.4){$_2$} \put(3.9,8.5){\small 3} \put(0.4,10.5){\small
4} \put(13,10.5){\small 5} \put(9.4,8.5){\small 6}
\put(7.85,0){\small 1}

\thicklines

\put(3,10){\line(-1,1){2}} \put(1,12){\line(1,0){12}}
\put(13,12){\line(-1,-1){2}} \put(3,10){\line(1,0){8}}

\thinlines

\put(3,10){\line(5,1){10}} \put(11,10){\line(-5,1){10}}
\put(7,0){\line(-1,2){6}} \put(7,0){\line(-2,5){4}}
\put(7,0){\line(2,5){4}} \put(7,0){\line(1,2){6}}

\put(4,0){\mbox{$S_1$}}


\put(21.77,9){$_{\bullet}$} \put(17.77,10){$_{\bullet}$}
\put(15.77,12){$_{\bullet}$} \put(27.77,12){$_{\bullet}$}
\put(25.77,10){$_{\bullet}$} \put(21.77,0){$_{\bullet}$}
\put(21.77,10.8){$_{\bullet}$}

\put(21.8,11.4){$_2$} \put(18.8,8.4){\small 3}
\put(15.4,10.5){\small 4} \put(28,10.5){\small 5}
\put(24.5,8.4){\small 6} \put(22.8,0){\small 1}
\put(21,7.8){\small 7}

\thicklines

\put(22,9){\line(-4,1){4}} \put(18,10){\line(-1,1){2}}
\put(16,12){\line(1,0){12}} \put(28,12){\line(-1,-1){2}}
\put(26,10){\line(-4,-1){4}}

\thinlines

\put(22,9){\line(0,1){1.8}} \put(18,10){\line(5,1){10}}
\put(26,10){\line(-5,1){10}}

\put(22,0){\line(-1,2){6}} \put(22,0){\line(-2,5){4}}
\put(22,0){\line(0,1){9}} \put(22,0){\line(2,5){4}}
\put(22,0){\line(1,2){6}}

\put(19,0){\mbox{$S_3$}}


\put(32.77,10){$_{\bullet}$} \put(30.77,12){$_{\bullet}$}
\put(42.77,12){$_{\bullet}$} \put(40.77,10){$_{\bullet}$}
\put(36.77,0){$_{\bullet}$} \put(36.77,10.8){$_{\bullet}$}

\put(36.8,11.4){$_2$} \put(33.8,9){$3$} \put(30.5,10.5){$4$}
\put(42.9,10.5){$5$} \put(39.5,9){$6$} \put(37.85,0){$1$}

\thicklines

\put(33,10){\line(-1,1){2}} \put(31,12){\line(1,0){12}}
\put(43,12){\line(-1,-1){2}}

\thinlines

\put(37,0){\line(0,1){11}}

\put(33,10){\line(5,1){10}} \put(41,10){\line(-5,1){10}}
\put(37,0){\line(-1,2){6}} \put(37,0){\line(-2,5){4}}
\put(37,0){\line(2,5){4}} \put(37,0){\line(1,2){6}}

\put(34,0){\mbox{$S_2$}}


\put(51.82,9){$_{\bullet}$} \put(47.82,10){$_{\bullet}$}
\put(45.92,11.9){$_{\bullet}$} \put(57.7,11.9){$_{\bullet}$}
\put(55.82,10){$_{\bullet}$} \put(51.825,0){$_{\bullet}$}

\put(51.7,9.7){$2$} \put(48.8,8.5){$3$} \put(45.5,10.5){$4$}
\put(57.9,10.5){$5$} \put(54.5,8.5){$6$} \put(52.85,0){$1$}

\thicklines

\put(52,9){\line(-4,1){4}} \put(48,10){\line(-1,1){2}}
\put(46,12){\line(1,0){12}} \put(58,12){\line(-1,-1){2}}
\put(56,10){\line(-4,-1){4}}

\thinlines

\put(52,9){\line(-2,1){6}} \put(52,9){\line(2,1){6}}
\put(52,0){\line(-1,2){6}} \put(52,0){\line(-2,5){4}}
\put(52,0){\line(0,1){9}} \put(52,0){\line(2,5){4}}
\put(52,0){\line(1,2){6}}

\put(44,7){\mbox{$=$}} \put(43.7,5){\mbox{$(\cong)$}}

\end{picture}

Observe that $S_2$ has two vertices of degree 5 but $S_1$ has no
degree 5 vertex. So, $S_1 \not \cong S_2$. Now, $S_3$ is obtained
from $S_1$ by starring at $7$ (in the edge $36$) and is obtained
from $S_2$ by starring at $7$ (in the edge $12$). Thus, $S_1
\approx S_2$ ($\approx S_3$).

\item[Bistellar Moves.]
Let $K$ be a $d$-dimensional ($d\geq 2$) pure simplicial complex
in $\RR^N$. Let $A$ be an $(d - k)$-simplex in $K$ such that ${\rm
lk}_K(A) = {\rm Bd}(B)$ for some $k$-simplex $B$ which is not in
$K$. Let $C = \langle A\cup B\rangle$. If $C \cap |K| = A \bdB$
then consider the simplicial complex $L= (K \setminus \{D : A <
D\in K\}) \cup \{FB : F<A, F\neq A\}$ (i.e., $L = (K\setminus
{\rm Cl}(A)\ast {\rm Bd}(B))\cup ({\rm Bd}(A)\ast {\rm Cl}(B)$).
We say that $L$ is obtained from $K$ by the {\em bistellar
$k$-move} $\kappa(A, B)$.

[For $0< k < d$, $C$ is a polytope of dimension $d$ or $d+1$ with
$d +2$ vertices. If $\dim(C) = d$ then $C= A \bdB =\bdA B$,
$\bdC=\bdA\bdB$ and $|L|=|K|$. If $\dim(C) = d+1$ then $C= A B$,
$\bdC=A\bdB \cup \bdA B$ and $|L|$ is pl homeomorphic to $|K|$.]

If $k = d$ we take $B=b\in\AAA$ then $\kappa (A, b)$ is
equivalent to starring at $b$. If $k = 0$ then $\kappa (u, B)$ is
equivalent to collapsing the vertex $u$. If $0 < k < d$ then
$\kappa(A, B)$ is called a {\em proper} bistellar move.

Observe that $S_2$ (defined above) is obtained from $S_1$ by the
bistellar 1-move $\kappa(36, 12)$.

\setlength{\unitlength}{3mm}

\begin{picture}(45,12)(0,-2.5)


\thicklines

\put(2,5){\line(4,3){4}} \put(2,5){\line(3,-1){3}}
\put(2,5){\line(1,-1){4}} \put(8,5){\line(-2,3){2}}
\put(8,5){\line(-3,-1){3}} \put(8,5){\line(-1,-2){2}}
\put(5,4){\line(1,4){1}} \put(5,4){\line(1,-3){1}}

\put(14,5){\line(4,3){4}} \put(14,5){\line(3,-1){3}}
\put(14,5){\line(1,-1){4}} \put(20,5){\line(-2,3){2}}
\put(20,5){\line(-3,-1){3}} \put(20,5){\line(-1,-2){2}}
\put(17,4){\line(1,4){1}} \put(17,4){\line(1,-3){1}}
\put(18,1){\line(0,1){3.1}} \put(18,4.6){\line(0,1){3.4}}

\put(9,6){\vector(1,0){4}} \put(13,4){\vector(-1,0){4}}
\put(32,6){\vector(1,0){4}} \put(36,4){\vector(-1,0){4}}

\put(24,4){\line(1,1){4}} \put(24,4){\line(2,-1){4}}
\put(28,2){\line(0,1){6}} \put(31,4){\line(-3,-2){3}}
\put(31,4){\line(-3,4){3}}

\put(37,4){\line(1,1){4}} \put(37,4){\line(2,-1){4}}
\put(41,2){\line(0,1){6}} \put(44,4){\line(-3,-2){3}}
\put(44,4){\line(-3,4){3}}

\thinlines

\put(2,5){\line(1,0){2.9}} \put(8,5){\line(-1,0){2.4}}

\put(14,5){\line(1,0){2.9}} \put(20,5){\line(-1,0){1.7}}

\put(17.5,4.7){\mbox{-}}

\put(24,4){\line(1,0){3.6}} \put(28.4,4){\line(1,0){2.6}}

\put(37,4){\line(1,0){3}} \put(41.4,4){\line(1,0){2.6}}

\put(40,5){\line(1,-3){1}} \put(40,5){\line(-3,-1){3}}
\put(40,5){\line(1,3){1}} \put(44,4){\line(-4,1){2.6}}

\put(40.2,4.8){\mbox{$.$}} \put(40.65,4.752){\mbox{$.$}}
\put(40.5,3.948){\mbox{$.$}}

\put(9,7){\mbox{{\small $1$-move}}} \put(9.3,2.5){\mbox{{\small
$2$-move}}}

\put(32,7){\mbox{{\small $0$-move}}} \put(32.3,2.5){\mbox{{\small
$3$-move}}}

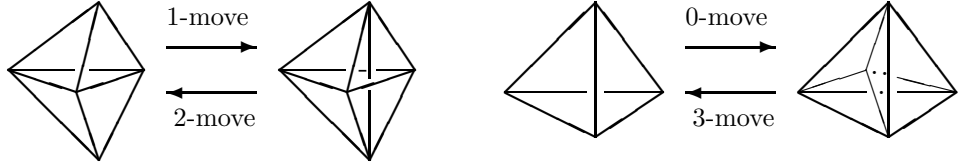
\put(14,-1.5){\mbox{Bistellar moves in dimension $3$}}

\end{picture}

Two pure simplicial complexes $M$ and $N$ are called {\em
bistellar equivalent} if there exists a finite sequence $M_1,
\dots, M_n$ of pure simplicial complexes such that $M_1=M$,
$M_n=N$ and $M_{i+1}$ is obtained from $M_{i}$ by a bistellar
move for $1\leq i\leq n-1$.

\item[Triangulations.]
A {\em triangulation} of a polyhedron $P$ is a pair $(K, t)$,
where $K$ is a simplicial complex and $t\colon |K|\to P$ is a pl
homeomorphism. Moreover, if $t$ is linear on each simplex then
$(K, t)$ is called a {\em linear triangulation}. We identify two
triangulations of a polyhedron $P$ if they differ by an
isomorphism (i.e., we identify $(K_1, t_1)$ and $(K_2, t_2)$ if
there is an isomorphism $i \colon K_1 \to K_2$ such that $t_2
\circ |i| =t_1$). If $(K, t)$ is a triangulation of $P$ and
$K^{\prime} \lhd K$ then $(K^{\prime}, t)$ is called a {\em
subdivision} of $(K, t)$.

For a simplicial complex $K$, if $|K|$ is homeomorphic to a
topological space $X$ then we say that $K$ is a {\em
triangulation} of $X$.

\item[Combinatorial Balls.]
For $d \geq 0$, let $\Delta^{d}$ be the $d$-simplex $\{(x_1, x_2,
\dots, x_{d + 1}) : x_i \geq 0 \mbox{ for } 1 \leq i \leq d + 1,
\sum_{i = 1}^{d + 1} x_i \leq 1\}$ in $\RR^{d + 1}$ with vertices
$v_1 = (1, 0, \dots, 0), \dots, v_{d + 1} = (0, \dots, 0, 1)$.
Then $\Delta^{d}$ is pl homeomorphic to the pl ball
$I^{\hspace{.2mm}d}$. Let ${\rm Cl}(\Delta^{d})$ denote the
simplicial complex whose simplices are all the faces of
$\Delta^{d}$. Then $|{\rm Cl} \Delta^{d}| = \Delta^{d}$. The
simplicial complex ${\rm Cl}(\Delta^{d})$ is called the {\em
standard $d$-ball}. A finite simplicial complex $K$ is called a
{\em combinatorial $d$-ball}, if $|K|$ is pl homeomorphic to
$\Delta^{d}$ (i.e., $K \approx {\rm Cl}(\Delta^{d})$).

\item[Combinatorial Spheres.]
Let ${\rm Bd}(\Delta^{d + 1})$ denote the simplicial complex ${\rm
Cl}(\Delta^{d + 1}) \setminus \{\Delta^{d + 1}\}$ ($d \geq 0$).
Then $|{\rm Bd}(\Delta^{d + 1})|$ is homeomorphic to the sphere
$S^{\,d}$. ($|{\rm Bd}(\Delta^{d + 1})| = \bdDelta$ is pl
homeomorphic to the pl sphere $\partial I^{d+1}$.) The simplicial
complex ${\rm Bd}(\Delta^{d + 1})$ is called the {\em standard
$d$-sphere} and is denoted by $S^{\,d}_{d + 2}(V)$ (or simply by
$S^{\,d}_{d + 2}$), where $V = \{v_1, \dots, v_{d + 2}\}$ is the
vertex-set of $\Delta^{d + 1}$.

A simplicial complex $K$ is called a {\em combinatorial
$d$-sphere}, if $|K|$ is pl homeomorphic to $|S^{\,d}_{d+2}|$
(i.e., by Proposition \ref{prop5}, $K \approx S^{\,d}_{d+2}$).

If a combinatorial $d$-sphere is $k$-neighbourly then, by
Corollary \ref{re4.6}, $k \leq \lfloor\frac{d+1}{2} \rfloor$.
Thus, a $\lfloor\frac{d+1}{2} \rfloor$-neighbourly combinatorial
$d$-sphere is called {\em neighbourly}.

\item[Polytopal Spheres.]
For $d\geq 0$, let $P$ be a simplicial $(d + 1)$-polytope in
$\RR^{d+1}$. Then the set of proper faces of $P$ form a
combinatorial $d$-sphere and is called the {\em boundary complex}
of the polytope $P$.

A combinatorial $d$-sphere $S$ is called a {\em polytopal sphere}
if it is isomorphic to the boundary complex of a simplicial $(d +
1)$-polytope.

\item[Stacked Spheres.] A combinatorial $d$-sphere $S$ is called
a {\em stacked} sphere if there is a sequence $S_1, \dots, S_k$ of
combinatorial $d$-spheres such that $S_1 = S^{\hspace{.4mm}d}_{d
+ 2}$, $S_k = S$ and $S_{j + 1}$ is obtained from $S_j$ by
starring a vertex on a facet of $S_j$ for $1 \leq j \leq k - 1$.

It follows from Proposition \ref{prop2} that a {\em stacked
$d$-sphere} is isomorphic to the boundary complex of a stacked $(d
+ 1)$-polytope. Clearly, the face-vector of an $n$-vertex stacked
$d$-sphere $S$ is given by
\begin{eqnarray} \label{stsphere}
f_k(S) & = & {d+2 \choose k+1} + (n-d-2) {d+1 \choose k} = {d+1
\choose k}n - {d+2 \choose k+1}k ~ \mbox{
for }  1\leq k < d, \nonumber \\
f_{d}(S) & = & (d+2) + (n - d - 2)d = dn - (d+2)(d-1).
\end{eqnarray}

\item[Combinatorial Manifolds.]
A simplicial complex $K$ is called a {\em combinatorial
$d$-manifold} if ${\rm lk}_K(v)$ is a combinatorial $(d-
1)$-sphere (i.e., $\approx S^{\,d-1}_{d +1}$) for each vertex $v$
in $K$. Clearly, a two dimensional complex $K$ is a combinatorial
2-manifold if the link of each vertex is a cycle. (A {\em cycle}
is a connected finite graph in which the degree of each vertex is
2. A cycle with $n$ vertices is called an $n$-cycle and is
denoted by $S^1_n$. An $n$-cycle with edges $v_1v_2, \dots,
v_{n-1}v_n$, $v_nv_1$ is also denoted by $S^1_n(v_1, \dots,
v_n)$.)

Since the link of a vertex in $S^{\hspace{.4mm}k}_{k + 2}$ is a
standard $(k - 1)$-sphere, it follows that $(i)$
$S^{\hspace{.4mm}k}_{k + 2}$ (and hence a combinatorial
$k$-sphere) is a combinatorial $k$-manifold and $(ii)$ if
$\sigma$ is an $i$-simplex in a combinatorial $d$-manifold $K$
then ${\rm lk}_K(\sigma)$ is a combinatorial $(d - i - 1)$-sphere
for $0 \leq i \leq d - 1$.

A simplicial complex $K$ is called a {\em combinatorial
$d$-manifold with boundary} if ${\rm lk}_K(v)$ is a combinatorial
$(d- 1)$-sphere or combinatorial $(d-1)$-ball for each vertex $v$
in $K$ and there exists a vertex $u$ whose link is a
combinatorial $(d-1)$-ball.

\item[Triangulated Manifolds.]
If the geometric carrier $|K|$ of a simplicial complex $K$ is a
closed topological $d$-manifold then $K$ is called a {\em
triangulated $d$-manifold}. So, a combinatorial manifold is
triangulated manifold and for $d\leq 3$, a triangulated
$d$-manifold is a combinatorial $d$-manifold.

\item[Homology Manifolds.] A $d$-dimensional simplicial complex
$K$ is called a {\em homology manifold} if for any $x \in |K|$ and
$i < d$, $H_i(|K|, |K|\setminus\{x\}; \ZZ) = 0$ and $H_d(|K|,
|K|\setminus\{x\}; \ZZ) = \ZZ$. So, a triangulated manifold is a
homology manifold.

\item[Eulerian Complexes.] A $d$-dimensional simplicial complex
$K$ is called an {\em Eulerian Complex} if $\chi({\rm
lk}_K(\sigma)) = 1 + (-1)^{d-i-1}$ for any $i$-simplex $\sigma$,
$0 \leq i< d$. So, a triangulation of a sphere is an Eulerian
Complex.

\item[PL Structures on Manifolds.]
A compact topological manifold $M$ is called {\em triangulable}
if it is homeomorphic to the geometric carrier of a simplicial
complex $K$. Moreover, if $K$ is a combinatorial manifold (i.e.,
by Proposition \ref{prop6}, $|K|$ is a pl manifold) then we say
$K$ is a {\em combinatorial triangulation} of $M$.

If a combinatorial manifold $K$ triangulates $M$, then the
combinatorial equivalence class ${\cal K}$ of combinatorial
manifolds containing $K$ is called a {\em combinatorial structure}
or {\em pl structure} of $M$.

\item[Pseudomanifolds.]
A pure $d$-dimensional simplicial complex $K$ is called a
$d$-dimensional {\em pseudomanifold} (or {\em $d$-pseudomanifold})
if (i) each $(d-1)$-simplex is a face of exactly two facets of $K$
and (ii) for any pair $\sigma_1$, $\sigma_2$ of facets of $K$,
there exists a sequence $\tau_1, \dots, \tau_n$ of facets of $X$,
such that $\tau_1 = \sigma_1$, $\tau_n = \sigma_2$ and $\tau_i
\cap \tau_{i+1}$ is a $(d-1)$-simplex of $K$ for $1 \leq i \leq
n-1$. By convention, $S^{\hspace{.2mm}0}_{2}$ is the only
0-pseudomanifold.

\item[Normal Pseudomanifolds.]
A $d$-pseudomanifold is said to be a {\em normal pseudomanifold}
if the links of all the simplices of dimension $\leq d - 2$ are
connected.  Clearly, the $1$-dimensional normal pseudomanifolds
are the cycles and the $2$-dimensional normal pseudomanifolds are
just the connected combinatorial $2$-manifolds. But, normal
pseudomanifolds of dimension $d$ form a broader class than
connected combinatorial $d$-manifolds for $d \geq 3$. In fact,
any connected triangulated manifold is a normal pseudomanifold.

\item[Irreducible Pseudomanifolds.] For $n \geq d + 3$, an $n$-vertex
$d$-pseudomanifold $M$ is called {\em irreducible} if $M$ can not
be written as $S^{\,c}_{c+2} \ast N$ for some pseudomanifold $N$
and $c\geq 0$. $M$ is called completely reducible if it is the
join of one or more standard spheres. By Theorem \ref{n=d+3},
$(d+3)$-vertex $d$-pseudomanifolds are completely reducible.

\item[One-Point Suspension.] Let $K$ be an $n$ vertex
$d$-dimensional pseudomanifold in $\RR^m \equiv
\RR^m\times\{0\}\subseteq \RR^{m+1}$ and $u \in V(K)$. Let $v =
(0, \dots, 0, 1) \in \RR^{m+1}$. Consider the $(d+1)$-dimensional
pseudomanifold $\Sigma_u K$ whose facet-set is $\{u\alpha :
\alpha$ a facet of $K$ and $u \not\in \alpha\}\cup \{v\beta :
\beta$ a facet of $K\}$. Observe that $uv$ is an edge of
$\Sigma_u K$ and if $w$ is an interior point in $uv$ then the
simplicial complex obtained from $\Sigma_u K$ by starring at $w$
is isomorphic to $K \ast S^0_2$. So, $|K|$ is homeomorphic to the
suspension of $|K|$. The pseudomanifold $\Sigma_u K$ is called
the {\em one-point suspension} of $K$ (see \cite{bd2} for more).

\item[Complementarity.] A simplicial complex $K$ is said to satisfy
{\em complementarity} if $\emptyset \neq U \subseteq V(K)$ and
$U\neq V(K)$ imply exactly one of $\langle U\rangle$, $\langle
V(K) \setminus U \rangle$ is a simplex of $K$. The simplicial
complexes $\RR P^2_6$ and $\CC P^2_9$ (in Examples \ref{eg2} and
\ref{eg10} respectively) satisfy complementarity.

\item[Abstract Simplicial Complex.] An {\em abstract simplicial
complex} is a collection of non-empty finite sets (sets of {\em
vertices}) such that every non-empty subset of a member is also a
member. For $i \geq 0$, a member of size $i+1$ is called an {\em
$i$-simplex} of the complex. For an abstract simplicial complex
$K$, $V(X)$ denotes the vertex-set of $K$. An abstract simplicial
complex is called {\em pure} if all the maximal simplices contain
same number of vertices. For an abstract simplicial complex $X$,
${\rm EG}(X)$ denote edge-graph of $X$ (i.e., ${\rm EG}(X)$
consists of vertices and edges of $X$). If ${\rm EG}(X)$ is
connected then we say that $X$ is {\em connected}.

If $K$ is a simplicial complex then $K_a := \{\sigma : \sigma
\subseteq V(K)$, $\langle\sigma\rangle \in K\}$ is an abstract
simplicial complex and is called the {\em abstract simplicial
complex corresponding to $K$}.

Let $X$ be a finite abstract simplicial complex. A simplicial
complex $K$ is called a {\em geometric realization} of $X$ (and
is denoted by $X_{gr}$) if $K_a$ is isomorphic to $X$. Clearly,
two geometric realizations of a finite abstract simplicial
complex are isomorphic.

Let $X$ be a finite abstract simplicial complex with $V(X) =
\{v_1, \dots, v_n\}$. Let $A$ be an $(n-1)$-simplex with vertices
$a_1, \dots, a_n$ in $\RR^{n-1}$. Let $K =\{a_{i_1}\cdots a_{i_k}
: \{v_{i_1}, \dots, v_{i_k}\}$ is a simplex of $X\}$. Then $K$ is
a subcomplex of the simplicial complex ${\rm Cl}(A)$. Clearly,
$K$ is a geometric realization of $X$. So, geometric realizations
exist for finite complexes. For infinite case see \cite{sp}.

\item[Isomorphism and Automorphism.] An {\em isomorphism} between
two abstract simplicial complexes $X$ and $Y$ is bijection
$\varphi \colon V(X) \to V(Y)$ such that $A \in X$ if and only if
$\varphi(A)\in Y$. Two abstract simplicial complexes are called
{\em isomorphic} if such an isomorphism exists. We identify two
isomorphic complexes. An isomorphism from an abstract simplicial
complex $X$ to itself is called an {\em automorphism}. All the
automorphisms of $X$ form a group under composition, which is
denoted by ${\rm Aut}(X)$.

\item[Quotient Complex and Proper Action.] Let $X$ be an abstract
simplicial complex and $G$ be a group of automorphism (i.e., $G$
is subgroup of ${\rm Aut}(X)$). Let $\eta : V(X) \to V(X)/G$ be
the natural projection. Let $X/G$ denote the abstract simplicial
complex $\{\eta(A) : A\in X\}$. This complex is called the {\em
quotient complex}.

Let $X$ be a connected abstract simplicial complex. For $u, v\in
V(X)$, let $d_X(u, v)$ denote the length of a shortest path from
$u$ to $v$ in ${\rm EG}(X)$. (Then $d_X$ is a metric on $V(X)$.)
Let $G$ be a group of automorphism of $X$. It is easy to see that
if $X$ is a pure $d$-dimensional abstract simplicial complex and
$d_X(u, g(u)) \geq 2$ for all $u \in V(X)$ and $1\neq g\in G$
then $X/G$ is also a pure $d$-dimensional abstract simplicial
complex. We say that $G$ acts {\em properly} on $X$ if $d_X(u,
g(u)) \geq 3$ for all $u\in V(X)$ and $1\neq g\in G$.

\end{description}

\begin{prop}$\!\!${\bf .} \label{prop1}
Let $K$ be a combinatorial $d$-manifold. Let $G$ be a group of
automorphism of $K_a$. Let $K/G$ denote the geometric realization
$(K_a/G)_{gr}$ of the quotient $K_a/G$. If $G$ acts properly on
$K_a$ then $K/G$ is also a combinatorial $d$-manifold.
\end{prop}

\noindent {\bf Proof.} Let $\eta : V(K_a) \to V(K_a/G)$ be the
natural projection. Then $\eta$ induces an abstract simplicial
map $\eta_{gr}$ from $K$ to $K/G$. For $v\in V(K_a)$, let $[v]=
\eta(v)$. Since the action is proper, it follows that $K/G$ is a
pure $d$-dimensional simplicial complex. Assume that $V(K_a) =
V(K)$ and $V(K/G) = V(K_a/G) = \{[v] \, : \, v\in K_a\}$.

Let $[u]$ be a vertex of $K/G$. Let $v$ and $w$ be two vertices in
${\rm lk}_K(u)$. Since $vu$ and $uw$ are edges in $K$, it follows
that the length of the shortest path in $K$ between $v$ and $w$
is at most 2. Therefore (since the action of $G$ is proper), $[v]
\neq [w]$. This implies that $\eta_{gr}|_{{\rm lk}_K(u)} : {\rm
lk}_K(u) \to {\rm lk}_{K/G}([u])$ is injective and hence an
isomorphism. Thus the link of each vertex in $K/G$ is a
combinatorial sphere. This proves the result. \hfill $\Box$

\begin{prop}$\!\!${\bf .} \label{prop2}
Let $M$ be the boundary complex of a simplicial $(d+1)$-polytope
$P$. Let $N$ be the combinatorial $d$-sphere obtained from $M$ by
starring a vertex in a facet $\sigma$ of $M$. Then $N$ is
isomorphic to the boundary of the polytope $Q$ which is obtained
from $P$ by attaching a $(d+1)$-simplex along the $d$-face
$\sigma$ of $P$.
\end{prop}

\noindent {\bf Proof.} Assume that $P$ is in $\RR^{\hspace{.2mm}d
+ 1}$. Let $\sigma= v_0\cdots v_{d+1}$ and $a = \frac{1}{d+2}(v_0
+ \cdots + v_{d+1})\in \sigmaint$. We may assume that $N$ is
obtained from $M$ by starring at $a$. Let $L$ be the closed half
line through $a$ and perpendicular to $\sigma$ such that $L \cap
P = \{a\}$. Then there exists $\varepsilon > 0$ such that $x \in L
\setminus \{a\}$ and distance between $a$ and $x$ $\leq
\varepsilon$ implies $P\cup (x\sigma)$ is convex. ($x\sigma$
denotes the join of $x$ and $\sigma$.) Fix a point $v$ on $L$ at
a distance $\varepsilon$ from $a$. Let $Q =P\cup (v\sigma)$. Then
$Q$ is a $(d+1)$-polytope. Let $\varphi$ be the map from $N$ to
the boundary of $Q$ given by $\varphi(u)= u$ if $u$ is a vertex of
$P$ and $\varphi(a) = v$. It is easy to see that $\varphi$ is an
isomorphism. \hfill $\Box$

\bigskip

In \cite{bd2}, we have shown the following\,:

\begin{prop}$\!\!${\bf .} \label{prop3}
Let $\Sigma_u K$ be the one-point suspension of a pseudomanifold
$K$. The pseudomanifold $\Sigma_u K$ is a polytopal sphere if and
only if $K$ is so.
\end{prop}

\bigskip

Here we present some basic results in pl-topology. See \cite{rs}
for proofs.

\begin{prop}$\!\!${\bf .} \label{prop4}
Any compact polyhedron is the geometric carrier of some
simplicial complex.
\end{prop}

\begin{prop}$\!\!${\bf .} \label{prop5}
Let $K$ and $L$ be two simplicial complexes. If $f\colon |K|\to
|L|$ is pl, then there are simplicial subdivisions $K^{\prime}
\lhd K$ and $L^{\prime}\lhd L$ such that $f\colon |K^{\prime}|\to
|L^{\prime}|$ is simplicial.
\end{prop}

\begin{prop}$\!\!${\bf .} \label{prop6}
Suppose $K$ is a simplicial complex then $|K|$ is a pl
$n$-manifold if and only if ${\rm lk}_K(v)\approx {\rm
Bd}(\Delta_n)$ or ${\rm Cl}(\Delta_{n-1})$ for each $v\in V(K)$.
\end{prop}

\begin{prop}$\!\!${\bf .} \label{prop7}
For $p$, $q\geq 1$, let $B^{\hspace{.2mm}p}$,
$B^{\hspace{.2mm}q}$ be combinatorial balls $($of dimensions $p$
and $q$ respectively$)$ and $S^{\hspace{.3mm}p-1}$,
$S^{\hspace{.3mm}q - 1}$ be combinatorial spheres $($of dimensions
$p - 1$ and $q - 1$ respectively$)$. Then {\rm (i)}
$B^{\hspace{.2mm}p} \ast B^{\hspace{.2mm}q}$ is a combinatorial
$(p + q + 1)$-ball, {\rm (ii)} $S^{\hspace{.3mm}p - 1} \ast
B^{\hspace{.2mm}q}$ is a combinatorial $(p + q)$-ball and {\rm
(iii)} $S^{\hspace{.3mm}p - 1} \ast S^{\hspace{.3mm}q - 1}$ is a
combinatorial $(p + q - 1)$-sphere.
\end{prop}


\section{A Brief History of Triangulations}

\begin{itemize}
\item
It was shown by Rado in 1924 that all 2-manifolds are
triangulable. Since the link of a vertex in a triangulated
2-manifold is a cycle, 2-manifolds have pl structures.
\item
In 1935, Cairns proved that each closed smooth manifold is
triangulable.
\item
In 1940, Whitehead proved that each closed smooth manifold has a
pl structure.
\item
In 1952, Moise showed that all 3-manifolds are triangulable.
Again, the link of a vertex in a triangulated 3-manifold is a
triangulation of the 2-sphere and all triangulations of the
2-sphere are combinatorial 2-spheres. So, 3-manifolds have pl
structures. Moise also showed that each 3-manifold admits a
unique pl structure.
\item
In 1960, Kervaire gave the example of a pl 10-manifold which is
not smoothable.
\item
In 1961, Eells and Kuiper (independently, Tamura) gave examples
of 8 dimensional pl manifolds which are not smoothable.
\item
In 1964, Lojaciewitz proved that each real algebraic variety is
triangulable.
\item
In 1967, Kuiper obtained algebraic equations for all
non-smoothable pl 8-manifolds.
\item
It is shown by Munkres in 1967 that there is an one to one
correspondence between the set of smooth homotopy $m$-spheres and
the set of pl homotopy $m$-spheres for $3\leq m \leq 4$. So, by
Freedman's classification of 4-manifolds, there is an one to one
correspondence between the set of pl structures on $S^{\,4}$ and
the set of smooth structures on $S^{\,4}$.
\item
In 1969, Kirby and Siebenmann (independently, Lashof and
Rosenberg) proved that (i) there is exactly one well defined
obstruction in $H^{\,4}(M; \ZZ_2)$ to imposing a pl structure on
a closed topological $m$-manifold $M$, $m \geq 5$ and (ii) given
one pl structure, there is a bijection between the class of
distinct pl structures and $H^{\,3}(M, \ZZ_2)$. Therefore,
$S^{\hspace{.2mm}m}$ has a unique pl structure for $m \geq 5$.
[For smooth structures on $S^{\hspace{.2mm}m}$ we know the
following: In 1963, Milnor and Kervaire proved that the set
$\Theta^{\rm DIFF}_m$ of smooth (homotopy) $m$-spheres is a
finite abelian group under the connected sum operation for $m
\geq 5$. For $m = 5, 6, 7, 8, 9, 10, 11$, $\Theta^{\rm DIFF}_m =
0, 0$, $\ZZ_{28}$, $\ZZ_2$, $\ZZ_2 \oplus \ZZ_2 \oplus\ZZ_2$,
$\ZZ_6$, $\ZZ_{992}$ respectively.]
\item
In 1970, Siebenmann showed that for each $n\geq 5$ there exists a
closed manifold $M^{\,n}$ of dimension $n$ which does not admit a
pl structure.
\item
In 1970, Siebenmann gave the example of a triangulable 5-manifold
which does not admit any pl structure.
\item
In 1974, Hirsch and Mazur showed that if the dimension of a closed
manifold $M$ with a pl structure is $\leq 7$ then $M$ is
smoothable. So, for $n\leq 7$, a $n$-dimensional closed
topological manifold is smoothable if and only if it has a pl
structure.
\item
In 1974, Hirsch showed that if $M\times N$ is smoothable, where
$M$ and $N$ are closed pl manifolds, then both $M$ and $N$ are
smoothable. So, if $M$ is an $8$-dimensional non-smoothable pl
manifold (by Eells and Kuiper such $M$ exists) then $M\times
S^{\,d}$ is a non-smoothable pl manifold of dimension $8+d$ for
all $d\geq 1$.
\item
In 1976, Galewski and Stern (independently, Matumoto) defined an
obstruction element $\tau\in H^5(M; \rho)$ such that the closed
topological $n$-manifold $M^{\,n}$, $n\geq 5$, is triangulable if
and only if $\tau=0$. Then they proved that all closed
topological manifolds of dimension $\geq 5$ are triangulable if
and only if there is a homology 3-sphere $\Sigma$ such that (i)
$\Sigma$ has Rohlin (or Rochlin) invariant 1 (i.e., bounding a
parallelizable 4-manifold of index 8) (ii)\footnote{This
condition is now superfluous because of Cannon's result mentioned
in Example \ref{eg18}} the $(n-3)$-fold suspension of $\Sigma$ is
homeomorphic to $S^{\,n}$ and (iii) $\Sigma\#\Sigma$ bounds a
smooth homology 4-disc. They also proved that  each simply
connected closed topological 6-manifold is triangulable.
\item
In 1977, Akbulut and King proved that each pl manifold of
dimension $\leq 10$ is homeomorphic to a real algebraic variety.
\item
In 1982, Freedman showed that there are closed 4-manifolds which
are not smoothable. So, by Hirsch and Mazur's result, there are
closed 4-manifolds which have no pl structures.
\item
In 1985, Casson showed that there exists a closed 4-manifold
which is not triangulable.
\item
For more, see \cite{ku}, \cite{si} and the following AMS
Mathematical Review numbers\,: 2,73e, 14,72d, 22\,\#12536,
25:\,\#2608, 25:\,\#2612, 26\,\#5584, 26:\,\#6978, 26\,\#6980,
31\,\#5209, 33\,\#6641, 35\,\#3671, 39\,\#3494, 39\,\#3500,
40\,\#895, 42\,\#6837, 54:\,\#3711, 54\,\#8650, 54\,\#11335,
55\,\#13434, 80e:\,57019, 80m:\,57014, 81b:\,57015, 81f:\,57012,
84e:\,57006.

\end{itemize}

\section{Examples}

\begin{itemize}
\item[$\ast$]
{\sf In this section, we present some combinatorial manifolds.
Most of these are vertex-minimal triangulations. We will discuss
about these in the next section. }
\item[$\ast$]
{\sf Since the facet-set of a pure simplicial complex determines
the simplicial complex, we identify a pure simplicial complex with
its facet-set in this section.}
\item[$\ast$] {\sf Whenever we say that $v_1 \cdots v_k$ is a simplex
then we mean that $v_1 \cdots v_k$ is the convex hull of $k$
affinely independent points $v_1, \dots, v_k$ in some $\RR^N$. }
\item[$\ast$] {\sf In the examples below, $v_1\cdots v_m$ and $u_1
\cdots u_n$ are two simplices in a simplicial complex $K$ and
$\{v_1, \dots, v_m\}\cap \{u_1, \dots, u_n\} = \emptyset$ mean we
have taken the vertices of $K$ in some $\RR^N$ such that
$v_1\cdots v_m \cap u_1\cdots u_n = \emptyset$. This is possible
since $K$ is finite. In fact, if $K$ has $r$ vertices then we can
consider $K$ in $\RR^{r-1}$ by considering the vertices of $K$ to
be affinely independent. }

\end{itemize}

\begin{eg}$\!\!${\bf .} \label{eg1}
{\rm For $c_1, \dots, c_n\geq 0$, $S^{\,c_1}_{c_1+2}\ast
\cdots\ast S^{\,c_n}_{c_n +2}$ is a combinatorial $(c_1+\cdots
+c_n+n-1)$-sphere on $c_1+\cdots +c_n+2n$ vertices.}
\end{eg}

\begin{eg}$\!\!${\bf .} \label{eg2}
{\rm Two combinatorial 2-manifolds of positive Euler
characteristics.
\begin{eqnarray*} {\cal I} &= & \{uu_iu_{i+1}, u_iu_{i+1}v_{i+3},
v_iv_{i+1}u_{i+3}, vv_iv_{i+1} :
       1\leq i\leq 5 \} \mbox{ and}   \\
\RR P^2_6 &= & \{ uu_iu_{i+1}, u_iu_{i+1}u_{i+3} : 1\leq i\leq
5\}.
\end{eqnarray*}
Additions in the subscripts are modulo 5. The geometric carrier
of ${\cal I}$ is the 2-sphere and it corresponds to the boundary
of the Platonic solid icosahedron. The geometric carrier of $\RR
P^2_6$ is the real projective plane. The complex $\RR P^2_6$ is
called the hemi-icosahedron. Observe that $\ZZ_2$ ($= \{ 1,
-1\}$) acts properly on the abstract simplicial complex ${\cal
I}_a$ by $(-1)u = v$, $(-1)v = u$, $(-1)u_i = v_i$, $(-1)v_i =
u_i$ and ${\cal I}/\ZZ_2 = \RR P^2_6$. }
\end{eg}

\begin{eg}$\!\!${\bf .} \label{eg3}
{\rm Two combinatorial 2-manifolds of Euler characteristic $0$.
\begin{eqnarray*}
{\cal T} &=& \{ w_{i}w_{i+1}w_{i+3}, w_{i}w_{i+2}w_{i+3} : 1\leq
i\leq 7\} \mbox{ and} \\
{\cal K} & = & \{u_1u_2v_1, u_1u_2v_2, u_1u_3v_1, u_1u_3v_3,
u_1u_4v_2, u_1u_4v_4, u_2u_3v_2, u_2u_3v_4,
\\ && \quad u_2u_4v_1, u_2u_4v_3, u_3u_4v_3, u_3u_4v_4, u_1v_3v_4,
u_2v_3v_4, u_3v_1v_2, u_4v_1v_2\}.
\end{eqnarray*}
Additions in the subscripts are modulo 7 in ${\cal T}$. The
geometric carrier of ${\cal T}$ is the torus and the geometric
carrier of ${\cal K}$ is the Klein bottle.}
\end{eg}

\begin{eg}$\!\!${\bf .} \label{eg4}
{\rm Two combinatorial 2-manifolds of negative Euler
characteristics.
\begin{eqnarray*}
M &= & \{ u_{1+p}u_{4+p}u_{7+p}, u_{i+3p}u_{j+3p}u_{k+3p} :
(i,j,k)\in
\{(1,2,5), (1,3,5), (1,3,4),   \\
&& ~~~~ (1,8,9), (1,6,8), (1,2,6), (2,3,6)\}, 0\leq p\leq 2\} \mbox{ and}   \\
N &= & \{uu_iu_{i+1}, u_iu_{i+1}u_{i+4}, u_iu_{i+2}u_{i+4},
u_iu_{i+3}u_{i+6} : 1\leq i\leq 9\}.
\end{eqnarray*}
Additions in the subscripts are modulo 9. The geometric carrier
of $M$ is the non-orientable surface of Euler characteristic $-3$
and the geometric carrier of $N$ is the non-orientable surface of
Euler characteristic $-5$. }
\end{eg}

\begin{eg}$\!\!${\rm {\bf :}} \label{eg5}
{\rm Five 8-vertex combinatorial $3$-spheres.
\begin{eqnarray*}
S^{\, 3}_{8, 35} & = & \{1234, 1267, 1256, 1245, 2345, 2356, 2367,
3467,
3456, 4567, \\
&& \quad 1238, 1278, 2378, 1348, 3478, 1458, 4578, 1568, 1678,
5678\}, \\
S^{\, 3}_{8, 36} & = & \{1234, 1256, 1245, 1567, 2345, 2356, 2367,
3467,
3456, 4567, \\
&& \quad 1268, 1678, 2678, 1238, 2378, 1348, 3478, 1458, 1578,
4578\}, \\
S^{\, 3}_{8, 37} & = & \{1234, 1256, 1245, 1457, 2345, 2356, 2367,
3467,
3456, 4567, \\
&& \quad 1568, 1578, 5678, 1268, 2678, 1238, 2378, 1348, 1478,
3478\}, \\
S^{\, 3}_{8, 38} & = & \{1234, 1237, 1267, 1347, 1567, 2345, 2367,
3467,
3456, 4567, \\
&& \quad 2358, 2368, 3568, 1268, 1568, 1248, 2458, 1478, 1578,
4578\} \mbox{ and } \\
S^{\, 3}_{8, 39} &= & \kappa(46, 357)(S^{\, 3}_{8, 4}).
\end{eqnarray*}
First four of these combinatorial manifolds are neighbourly and
were found by Gr\"{u}nbaum and Sreedharan (in \cite{gs}, these
are denoted by $P^8_{35}$, $P^8_{36}$, $P^8_{37}$ and ${\cal M}$
respectively). They showed that $S^{\, 3}_{8, 35}$, $S^{\, 3}_{8,
36}$, $S^{\, 3}_{8, 37}$ are polytopal spheres and $S^{\, 3}_{8,
38}$ is a non-polytopal sphere (known as the {\bf
Br\"{u}ckner-Gr\"{u}nbaum sphere}). The sphere $S^3_{8, 39}$
(obtained from $S^{\, 3}_{8, 38}$ by the bistellar $2$-move
$\kappa(46, 357)$) is a non-polytopal sphere and found by
Branette in \cite{b3}.}
\end{eg}

\begin{eg}$\!\!${\bf .} \label{eg6}
{\rm Consider the 11-vertex pure 3-dimensional simplicial complex
$\RR P^{\,3}_{11}$ (on the vertex-set $\{1, \dots, 9, a, b\}$)
whose maximal simplices are
\begin{eqnarray*}
&& 1237, 123b, 1269, 126b, 1279, 135a, 135b, 137a, 1479, 147a,
1489, 148a,
1568, 156b, \\
&& 158a, 1689, 2348, 234b, 2378, 246a, 246b, 248a, 2578, 2579,
258a, 259a,
269a, 3459, \\
&& 345b, 3489, 359a, 3678, 367a, 3689, 369a, 4567, 456b, 4579,
467a, 5678.
\end{eqnarray*}
This simplicial complex  is a combinatorial 3-manifold and
triangulates the 3-dimensional real projective space $\RR P^3$.
This was first constructed by Walkup in \cite{w}. Theorem
\ref{walkup2} shows that 11 is the minimal number of vertices
required to triangulate $\RR P^{\hspace{.2mm}3}$.}
\end{eg}

\begin{eg}$\!\!${\bf .} \label{eg7}
{\rm Let $L^3_{12}$ be the 12-vertex pure 3-dimensional simplicial
complex (on the vertex set $\{1, \dots, 9, a, b, c\}$) whose
facets are
\begin{eqnarray*}
&& 1234, 123a, 1249, 1256, 1259, 126b, 12ab, 1347, 1378,
138a, 1479, 156c, 1579, 157c, \\
&& 16bc, 178c, 18ab, 18bc,
234c, 23ac, 2489, 248c, 2568, 2589, 2678, 267b, 278c, 27ab, \\
&& 27ac, 3456, 345b, 3467, 34bc, 3568, 3589, 359b,
3678, 389a, 39ac, 39bc, 456a, 45ab,\\
&&  4679, 469a, 489a, 48ab, 48bc, 56ac, 579b, 57ab, 57ac, 679b,
69ac, 69bc.
\end{eqnarray*}
This complex  is a combinatorial 3-manifold and triangulates the
lens space $L(3, 1)$ (\cite{lu1}). Since $L(3, 1)$ is a
$\ZZ_2$-homology $3$-sphere ($H_1(L(3,1), \ZZ)=\ZZ_3$,
$H_2(L(3,1), \ZZ)= 0$), it follows from Theorem \ref{z2hsa} that
12 is the least number of vertices required to triangulate $L(3,
1)$.}
\end{eg}

\begin{eg}$\!\!${\bf .} \label{eg8}
{\rm Consider the 15-vertex 3-dimensional pure simplicial complex
$$
T^3_{15} = \{u_i u_{i+p} u_{i+p+q} u_{i+p+q+r} : \{p, q, r\} =
\{1, 2, 4\}, 1 \leq i \leq 15\}.
$$
(Additions in the subscripts are modulo 15.) This simplicial
complex is a combinatorial 3-manifold and triangulates
$S^{\,1}\times S^{\,1}\times S^{\,1}$ (\cite{kl2}). A
generalization of this is presented in Example \ref{eg21}.}
\end{eg}

\begin{eg}$\!\!${\bf .} \label{eg9}
{\rm Let $H^3_{16}$ be the 16-vertex pure 3-dimensional
simplicial complex (on the vertex set $\{1, \dots, 9, a, b, c, d,
e, f, g\}$) whose facets are
\begin{eqnarray*}
&1249, 124f, 126e, 126f, 129e, 134c, 134f, 137a, 137c, 13af,
149c, 156d, 156e, 158b, 158d,& \\
& 15be, 16df, 178a, 178b, 17bc, 18ad, 19bc, 19be, 1adf, 235a,
235b, 237a, 237d, 23bd, 249d, & \\
& 24bd, 24bf, 258b, 258c, 25ac, 26ac, 26ae, 26cf, 279d, 279e,
27ae, 28bf, 28cf, 345e, 345f, & \\
& 34ce, 35af, 35be, 37cd, 3bde, 3cde, 4567, 456e, 457f, 467b,
46ab, 46ae, 47bf, 489c, 489d, & \\
& 48ad, 48ae, 48ce, 4abd, 567d, 579d, 579f, 589c, 589d, 59ac,
59af, 67bc, 67cd, 6abc, 6cdf,  & \\
& 78ae, 78bf, 78ef, 79ef, 8cef, 9abc, 9abg, 9afg, 9beg, 9efg,
abdg, adfg, bdeg, cdef, defg.&
\end{eqnarray*}
This simplicial complex is a combinatorial 3-manifold and was
constructed by Bj\"{o}rner and Lutz in \cite{bl}. The complex
$H^{\hspace{.2mm}3}_{16}$ has $f$-vector $(16, 106, 180, 90)$ and
triangulates the Poincar\'{e} homology 3-sphere. It follows from
Theorem \ref{z2hsa} that at least 12  vertices are required to
triangulate the Poincar\'{e} homology 3-sphere. }
\end{eg}

\begin{eg}$\!\!${\bf .} \label{eg10}
{\rm Consider a 9-vertex abstract simplicial complex ${\cal X}$
as follows. The vertices of $\cal X$ are the points of the affine
plane $\cal P$ over the 3-element field. Fix a set $\Pi$ of three
mutually parallel lines of $\cal P$ (i.e., $\Pi$ is a parallel
class of lines in $\cal P$). Let the lines in $\Pi$ be
$\gamma_0$, $\gamma_1$, $\gamma_2$, in a fixed cyclic
orientation. The set of maximal simplices of $\cal X$ is as
follows.
$$
\{\gamma_{i+1} \cup \gamma_{i} \setminus \{x\} \, : \, x \in
\gamma_i, 0\leq i\leq 2\} \cup \{\alpha \cup \beta \, : \, \alpha
\neq \beta \mbox{ two intersecting lines of } {\cal P} \mbox{
outside } \Pi\}.
$$
(Addition in the suffix is modulo 3.) This gives $3 \times 3 +
(9\times 6)/2 = 9 + 27= 36$ maximal simplices of ${\cal X}$. Then
the geometric realization ${\cal X}_{gr}$ of ${\cal X}$ is a
combinatorial 4-manifold. This ${\cal X}_{gr}$ triangulates the
complex projective plane and is denoted by $\CC P^2_9$
(\cite{bd1, k1, k2, my}). This was first constructed by
K\"{u}hnel and Banchoff in \cite{kb}. Check that, the link of any
vertex in $\CC P^2_9$ is isomorphic the Br\"{u}ckner-Gr\"{u}nbaum
sphere $S^{\,3}_{8, 38}$. }
\end{eg}

\begin{eg}$\!\!${\bf .} \label{eg11}
{\rm Consider the 11-vertex 4-dimensional pure simplicial complex
$S^{\,2,2}_{11}$ whose facets are
\begin{eqnarray*}
&&12346, 12347, 12369, 12379, 12458, 12459, 12468, 12479,
12568, 12569, 13467, 13567, \\
&& 13569, 1357a, 1359b, 135ab, 1379a, 139ab, 1458a, 1459b, 145ab,
1467b,
1468a, 146ab,\\
&& 1479b, 15678, 1578a, 1678b, 168ab, 178ab, 179ab, 23468, 23478,
2357a,
2357b, 235ab, \\
&& 2368a, 2369a, 2378b, 2379a, 238ab, 24589, 24789, 2568b, 2569a,
256ab,
25789, 2578b, \\
&& 2579a, 268ab, 3467b, 3468a, 3469a, 3469b, 3478b, 3489a, 3489b,
3567b,
3569b, 389ab, \\
&& 4569a, 4569b, 456ab, 4589a, 4789b, 5678b, 5789a, 789ab.
\end{eqnarray*}
The simplicial complex $S^{\,2,2}_{11}$  is a combinatorial
4-manifold and triangulates $S^{\,2}\times S^{\,2}$. This was
constructed by Lutz in \cite{lu1}. Observe that, by Theorems
\ref{d=4a} and \ref{d=4b}, 11 is the minimum number of vertices
required to triangulate $S^{\,2}\times S^{\,2}$. }
\end{eg}

\begin{eg}$\!\!${\bf .} \label{eg12}
{\rm Let $G$ be the subgroup of $S_{16}$ generated by $(2, 7)(4,
10)(5, 6)(11, 12)$ and $(1, 2, 3, 4, 5, 10)(6, 8, 9)(11, 12, 13,
14, 15, 16)$. Then $G \cong S_6$ and $G$ acts on the set $V=
\{u_1, \dots, u_{16}\}$ by $\alpha(u_{i}) = u_{\alpha(i)}$ for
$\alpha\in G$. This action induces an action on the set of
subsets of $V$, namely, $\alpha(U) = \{\alpha(a) : a\in U\}$ for
$U\subseteq V$. Consider the 16-vertex abstract simplicial complex
$$
K = \{\alpha(\{u_{1}, u_{2}, u_{4}, u_{5}, u_{11}\}), ~
\alpha(\{u_{1}, u_{2}, u_{4}, u_{11}, u_{13}\}) : \alpha\in G\}.
$$
Let $\RR^{\,4}_{16}= K_{gr}$ be the geometric realization of $K$.
Then $\RR^{\,4}_{16}$ is a combinatorial 4-manifold with
$f$-vector $(16, 120, 330, 375, 150)$ and triangulates the
4-dimensional real projective space $\RR P^{\,4}$. This was
constructed by Lutz in \cite{lu1}. It follows, from Theorem
\ref{rpdn}, that 16 is the minimum number of vertices required to
triangulate $\RR P^{\,4}$.}
\end{eg}

\begin{eg}$\!\!${\bf .} \label{eg13}
{\rm Let $\FF_4 = \{0, 1, x, y\}$ be the field of order 4.
Consider the space $\FF_4 \oplus \FF_4 = $ {\small $\left\{v_{0}
= \left(\hspace{-1.6mm}\begin{tabular}{c}
  0 \\
  0
\end{tabular}\hspace{-1.6mm}\right)\hspace{-1mm},
  \, v_{1} =
\left(\hspace{-1.6mm}\begin{tabular}{c}
  1 \\
  0
\end{tabular}\hspace{-1.6mm}\right)\hspace{-1mm},
  \,  v_{2} =
\left(\hspace{-1.6mm}\begin{tabular}{c}
  x \\
  0
\end{tabular}\hspace{-1.6mm}\right)\hspace{-1mm},
  \,   v_{3} =
\left(\hspace{-1.6mm}\begin{tabular}{c}
  y \\
  0
\end{tabular}\hspace{-1.6mm}\right)\hspace{-1mm},
  \,  v_{4} =
\left(\hspace{-1.6mm}\begin{tabular}{c}
  0 \\
  1
\end{tabular}\hspace{-1.6mm}\right)\hspace{-1mm},
 \,  v_{5} = \left(\hspace{-1.6mm}\begin{tabular}{c}
  1 \\
  1
\end{tabular}\hspace{-1.6mm}\right)\hspace{-1mm},
  ~  v_{6} =
\left(\hspace{-1.6mm}\begin{tabular}{c}
  x \\
  1
\end{tabular}\hspace{-1.6mm}\right)\hspace{-1mm},
  \,  v_{7} =
\left(\hspace{-1.6mm}\begin{tabular}{c}
  y \\
  1
\end{tabular}\hspace{-1.6mm}\right)\hspace{-1mm},\right.$
\newline
 $\left.  v_{8} = \left(\hspace{-1.6mm}\begin{tabular}{c}
  0 \\
  x
\end{tabular}\hspace{-1.6mm}\right)\hspace{-1mm},
  \,  v_{9} =
\left(\hspace{-1.6mm}\begin{tabular}{c}
  1 \\
  x
\end{tabular}\hspace{-1.6mm}\right)\hspace{-1mm},
  \,  v_{10} =
\left(\hspace{-1.6mm}\begin{tabular}{c}
  x \\
  x
\end{tabular}\hspace{-1.6mm}\right)\hspace{-1mm},
  \,  v_{11} =
\left(\hspace{-1.6mm}\begin{tabular}{c}
  y \\
  x
\end{tabular}\hspace{-1.6mm}\right)\hspace{-1mm},
  \,  v_{12} =
\left(\hspace{-1.6mm}\begin{tabular}{c}
  0 \\
  y
\end{tabular}\hspace{-1.6mm}\right)\hspace{-1mm},
 \, v_{13} = \left(\hspace{-1.6mm}\begin{tabular}{c}
  1 \\
  y
\end{tabular}\hspace{-1.6mm}\right)\hspace{-1mm},
  \,  v_{14} =
\left(\hspace{-1.6mm}\begin{tabular}{c}
  x \\
  y
\end{tabular}\hspace{-1.6mm}\right)\hspace{-1mm},
  \,  v_{15} =
   \left(\hspace{-1.6mm}\begin{tabular}{c}
  y \\
  y
\end{tabular}\hspace{-1.6mm}\right)\right\}$.} \newline
Let $G$ be the group generated by all the translations in $\FF_4
\oplus \FF_4$ and the matrix $A =$ {\small
$\left(\begin{array}{cc}
  0 & y \\
  y & 1
\end{array}\right)$}. Then the order of $G$ is 240 and $G$ acts
transitively on $\FF_4 \oplus \FF_4$. This action induces an
action on the set of subsets of $\FF_4 \oplus \FF_4$, namely,
$g(U) := \{g(a) : a\in U\}$ for $U \subseteq \FF_4 \oplus \FF_4$.

Consider the abstract simplicial complex ${\cal K}$ on the
vertex-set $\FF_4 \oplus \FF_4$ as
$$
{\cal K} = \{g(\{v_1, v_2, v_3, v_4, v_8\}), ~ g(\{v_1, v_4, v_6,
v_9, v_{10}\}) ~ : ~ g\in G\}.
$$
(One orbit of 4-simplices of length 240 and one orbit of
4-simplices of length 48.) Let ${\rm K3}_{16} = {\cal K}_{gr}$ be
the geometric realization of ${\cal K}$. Then ${\rm K3}_{16}$ is
a combinatorial 4-manifold and triangulates a K3 surface. This
was constructed by Casella and K\"{u}hnel in \cite{ck}. Since the
Euler characteristic of a K3 surface is 24, by Theorem
\ref{d=4a}, 16 is the minimum number of vertices required to
triangulate a K3 surface. }
\end{eg}

\begin{eg}$\!\!${\bf .} \label{eg14}
{\rm Consider the 12-vertex 5-dimensional pure simplicial complex
$S^{\,3,2}_{12}$ (on the vertex set $\{1, \dots, 9, a, b, c\}$)
whose facets are
\begin{eqnarray*}
&12346a, 12346b, 123478, 12347b, 12348a, 12357b, 12357c, 12359b,
12359c, 1236ab, 12378c,  & \\
& 1238ac, 1239ab, 1239ac, 124678, 12467b, 124689, 12469a, 12489a,
1257bc, 1259bc, 12678c, & \\
& 1267bc, 12689c, 1269ab, 1269bc, 1289ac, 134678, 13467b, 13468a,
13579b, 13579c, 13678c, & \\
& 1367bc, 1368ac, 136abc, 1379ab, 1379ac, 137abc, 145689, 14568a,
14569a, 14589c, 1458ac, & \\
& 1459ac, 1489ac, 15689b, 1568ab, 1569ab, 1579ab, 1579ac, 157abc,
1589bc, 158abc, 1689bc, & \\
& 168abc, 23456a, 23456c, 23458a, 23458b, 2345bc, 2346bc, 23478b,
235678, 23567c, 23568a, & \\
&23578b, 2359bc, 23678c, 2368ac, 236abc, 239abc, 24567a, 24567c,
24578a, 24578b, 2457bc, & \\
& 246789, 24679a, 2467bc, 24789a, 25678a, 26789a, 2689ac, 269abc,
345679, 34567c, 345689, & \\
& 34568a, 34579c, 34589b, 3459bc, 346789, 3467bc, 34789b, 3479bc,
356789, 35789b, 379abc, & \\
& 45679a, 4578ab, 4579ac, 457abc, 4589bc, 458abc, 4789ab, 479abc,
489abc, 56789a, 5689ab, & \\
& 5789ab, 689abc. &
\end{eqnarray*}
The simplicial complex $S^{\,3,2}_{12}$  is a combinatorial
5-manifold and triangulates $S^{\,3}\times S^{\,2}$. This was
constructed by Lutz in \cite{lu1}. Observe that, by Theorem
\ref{n<2d+4+i}, 12 is the minimum number of vertices required to
triangulate $S^{\,3}\times S^{\,2}$. }
\end{eg}

\begin{eg}$\!\!${\bf .} \label{eg15}
{\rm In $\RR^{\,d+1}$ consider the moment curve $M_{d+1}$ defined
parametrically by $x(t) = (t, t^2, \dots, t^{d+1})$. Let $t_1 <
t_2 < \cdots < t_{n}$ and $v_i = x(t_i)$ for $1\leq i\leq n$. For
$n\geq d+2$, let $V =\{v_1, \dots, v_{n}\}$. Let $C(n, d+1)$ be
the convex hull of $V$. Then $C(n, d+1)$ is a simplicial convex
$(d+1)$-polytope. The boundary complex of $C(n, d+1)$ is called
the {\em cyclic $d$-sphere} and is denoted by $C_n^{d}$. Then (i)
$C(n, d+1)$ (and hence $C_n^{d}$) is $\lfloor \frac{d+1}{2}
\rfloor$-neighbourly and (ii) a set $U$ ($\subseteq V$) of $d+1$
vertices spans a $d$-face of $C(n, d+1)$ if and only if any two
points of $V \setminus U$ are separated on $M_{d+1}$ by even
number of points of $U$ (see \cite[Pages 61--63]{gr}). Observe
that the link of a vertex in $C^{\,2c+1}_{m+1}$ is isomorphic to
$C^{\,2c}_{m}$.

If $d$ is odd then, by (ii), $v_1 v_3 \dots v_{d + 2}$ is not a
simplex of $C^{\,d}_n$. So, $C^{\,d}_n$ is not $(\frac{d + 1}{2} +
1)$-neighbourly. If $d$ is even then, by (ii), $v_2v_4 \dots v_{d
+ 2}$ is not a simplex of $C^{\,d}_n$. So, $C^{\,d}_n$ is not
$(\frac{d}{2} + 1)$-neighbourly. Thus, $C^{\,d}_n$ is not
$\lfloor\frac{d + 1}{2} + 1 \rfloor$-neighbourly for all $d\geq
1$.

For odd $d \geq 1$, consider the following  pure abstract
simplicial complex ${\cal C}$. The vertices of ${\cal C}$ are the
vertices of the $n$-vertex circle $S^{\, 1}_{n}$ and a set of
$d+1$ vertices is a maximal simplex of ${\cal C}$ if and only if
the induced subgraph of $S^{\, 1}_{n}$ on these $d+1$ vertices
has no connected component of odd size. By (ii), it follows that
$C^{\,d}_n$ is the geometric realization of this ${\cal C}$. }
\end{eg}

\begin{eg}$\!\!${\bf .} \label{eg16}
{\rm For $d \geq 2$ and $n \geq 2d + 3$, let $v_1, \dots, v_n$ be
$n$ affinely independent points (in $\RR^{\,n-1}$). Consider the
$(d+1)$-dimensional pure simplicial complex $X$ on the vertex set
$\{v_1, \dots, v_n\}$ given by\,:
$$
X = \{v_iv_{i+1}\cdots v_{i+d+1} \, : \, 1\leq i \leq n\}.
$$
(Addition in the suffix is modulo $n$.) Then $|X|$ is a pl
manifold with boundary. Let $K^{\hspace{.2mm}d}_{n}$ be the
boundary complex of $|X|$. More explicitly, the facet-set of
$K^{\hspace{.2mm}d}_{n}$ is
$$
\{\alpha : \alpha \mbox{ is a $d$-simplex in $X$ and $\alpha$ is
in a unique facet of } X\}.
$$
Then $K^{\hspace{.2mm}d}_{n}$ is a combinatorial $d$-manifold. It
was shown in \cite{k1} the following\,: (i)
$K^{\hspace{.2mm}d}_{2d + 3}$ triangulates $S^{\hspace{.2mm}d -
1} \times S^{\hspace{.2mm}1}$ for $d$ even, and triangulates the
twisted product $\TPSSD$ (the twisted $S^{\hspace{.2mm}d -
1}$-bundle over $S^{\hspace{.2mm}1}$) for $d$ odd. (ii)
$K^{\hspace{.2mm}d}_{2d+4}$ triangulates $S^{\hspace{.2mm}d-1}
\times S^{\hspace{.2mm}1}$ for all $d$. In particular,
$K^{\hspace{.2mm}3}_9$ triangulates the twisted product $\TPSS$
(often called the {\em $3$-dimensional Klein bottle}) and
$K^{\hspace{.2mm}3}_{10}$ triangulates the product
$S^{\hspace{.2mm}2} \times S^{\hspace{.2mm}1}$. The combinatorial
3-manifolds $K^{\hspace{.2mm}3}_9$ and $K^{\hspace{.2mm}3}_{10}$
were first constructed by Walkup in \cite{w}.

From the definition, it follows that the abstract simplicial
complex $(K^{\hspace{.2mm}d}_{n})_{a}$ corresponding to
$K^{\hspace{.2mm}d}_{n}$ is the pure abstract simplicial complex
whose vertices are the vertices of the $n$-cycle
$S^{\hspace{.2mm}1}_{n}(v_1, \dots, v_n)$ and the $d$-simplices
are the sets of $d+1$ vertices obtained by deleting an interior
vertex from the $(d + 2)$-paths in the cycle. (In fact, the
maximal simplices of $X_a$ are the $(d + 2)$-paths in
$S^{\hspace{.2mm}1}_{n}(v_1, \dots, v_n)$.) }
\end{eg}

\begin{eg}$\!\!${\bf .} \label{eg17}
{\rm For $d \geq 2$, consider the $(3d+5)$-vertex abstract
simplicial complexes ${\cal M}$ and ${\cal N}$ on the vertex set
$V = \{1, \dots, 3d + 5\}$ given by\,:
\begin{eqnarray*}
{\cal N} & = & \{\{i, \dots, j-1, j+1, i+d+1 \, : \, i+1 \leq j
\leq i+d, \, 1\leq i \leq 2d+4\}, \\
{\cal M} & = & {\cal N} \cup \{\{1, \dots, d+1\}, \{2d+5, \dots,
3d+5\}\}.
\end{eqnarray*}

Then ${\cal M}_{gr}$ is a stacked $d$-sphere and (hence) ${\cal
N}_{gr}$ triangulates $S^{\hspace{.2mm}d-1}\times [0, 1]$.

Let $p = (p_1,\dots, p_k)$ be a partition of $d + 1$. Put $s_0 =
0$ and $s_j = \sum_{i=1}^{j} p_i$ for $1 \leq j \leq k$. Let
$\pi_p$ be the permutation of $\{1, 2, \dots, d+1\}$ which is the
product of $k$ disjoint cycles $(s_{j - 1} + 1, s_{j - 1} +2,
\dots, s_{j})$, $1\leq j\leq k$. Since $\pi_p(i) \leq i + 1$ for
$1 \leq i \leq d + 1$, it follows that $d_{\cal M}(2d + 4 + i,
\pi_p(i)) \geq 3$. Consider the equivalence relation $\rho(p)$ on
$V$ given by\,: $j \rho(p) j$ for $j \in V$ and $(2d + 4 + i)
\rho(p) \pi_p(i)$ for $1 \leq i \leq d + 1$. Let $\eta_p \colon V
\to V/\rho(p)$ be the canonical surjection.

Let ${\cal K}(p) = {\cal N}/\rho(p)$ denote the abstract
simplicial complex whose vertex set is $V/\rho(p)$ and simplices
are $\eta_p(\sigma)$, where $\sigma \in {\cal N}$. Let $K^{d}_{2d
+ 4}(p) = {\cal K}(p)_{gr}$ be the geometric realization of
${\cal K}(p)$. It was shown in \cite{bd7} that $K^{d}_{2d +
4}(p)$ is a combinatorial $d$-manifold and triangulates
$S^{\hspace{.2mm}d - 1} \times S^{\hspace{.2mm}1}$ (respectively,
the twisted product $\TPSSD$) if $p$ is an even (respectively,
odd) partition of $d + 1$. If $p_0 = (1, 1, \dots, 1)$ then
$\pi_{p_0} = {\rm Id}$ and the corresponding combinatorial
$d$-manifold $K^{d}_{2d + 4}(p_0)$ is same as $K^{d}_{2d + 4}$
defined in Example \ref{eg16}.

[Recall that for any positive integer $n$, a {\em partition} of
$n$ is a finite weakly increasing sequence of positive integers
adding to $n$. The terms of the sequence are called the {\em
parts} of the partition. A partition of $n$ is {\em even}
(respectively, {\em odd}) if it has an even (respectively, odd)
number of even parts. It was shown in \cite{bd7} that the number
of odd (respectively, even) permutations of $n$ is $\geq
\frac{1}{2}\times$(the number of partitions of $n-1$).]}
\end{eg}

\begin{eg}$\!\!${\bf .} \label{eg18}
{\rm Let $H^3_{16}$ be the combinatorial 3-manifold defined in
Example \ref{eg9} and $u$ be a vertex of $H^3_{16}$. Let
$\Sigma^1H^3_{16} := \Sigma_uH^1_{16}$ be the one-point
suspension of $H^3_{16}$ and for $n\geq 2$, let $\Sigma^nH^3_{16}
:= \Sigma_u(\Sigma^{n-1}H^3_{16})$. Let $U = V(\Sigma^nH^3_{16})
\setminus V(H^3_{16})$. Then $\sigma := \langle U\rangle$ is an
$(n-1)$-simplex of $\Sigma^nH^3_{16}$ and ${\rm lk}_{\Sigma^n
H^3_{16}}(\sigma) = H^3_{16}$. Thus, $\Sigma^nH^3_{16}$ is not a
combinatorial manifold. Since $|H^3_{16}|$ is a homology 3-sphere
and $|\Sigma^nH^3_{16}|$ is the $n$-th suspension of $|M|$, by
Cannon's theorem (which states that the double suspension of any
homology $d$-sphere is homeomorphic to $S^{\,d + 2}$) \cite{ca},
$|\Sigma^nH^3_{16}|$ is homeomorphic to $S^{\,3 + n}$ for $n \geq
2$. So, $\Sigma^{n}H^{3}_{16}$ is a triangulation of $S^{\,3 + n}$
for $n \geq 2$. Clearly, $\Sigma^nH^3_{16}$ has $3 + n$ vertices.
So, for $d \geq 5$, $S^{\,d}$ has a $(d + 13)$-vertex
non-combinatorial (non-pl) triangulation. }
\end{eg}

\begin{eg}$\!\!${\bf .} \label{eg19}
{\rm Let $S^{\,3}_{8, 38}$ be the Br\"{u}ckner-Gr\"{u}nbaum
3-sphere defined in Example \ref{eg5} and $u$ be a vertex of
$S^{\,3}_{8, 38}$. Let $\Sigma^1S^{\,3}_{8, 38} := \Sigma_u
S^{\,3}_{8, 38}$ be the one-point suspension of $S^{\,3}_{8, 38}$
and for $n \geq 2$, let $\Sigma^n S^{\,3}_{8, 38} :=
\Sigma_u(\Sigma^{n - 1}S^{\,3}_{8, 38})$. Then
$\Sigma^nS^{\,3}_{8, 38}$ is an $(n + 8)$-vertex combinatorial
$(n+3)$-sphere. Since $S^{\,3}_{8, 38}$ is a non-polytopal sphere,
by Proposition \ref{prop3}, $\Sigma^n S^{\,3}_{8, 38}$ is a
non-polytopal sphere. So, for $d \geq 3$, there exists an $(d +
5)$-vertex non-polytopal combinatorial $d$-sphere. Note that, by
Theorem \ref{n=d+4b}, a $(d + k)$-vertex combinatorial $d$-sphere
is a polytopal sphere for $2 \leq k \leq 4$.}
\end{eg}

\begin{eg}$\!\!${\bf .} \label{eg20}
{\rm For $d\geq 1$, let $\Delta^{d+1}$ be the $d$-simplex with
vertices $v_1= (1, 0. \dots, 0)- \frac{1}{d+2}(1, 1, \dots, 1)$,
$\dots, v_{d+2} = (0, \dots, 0, 1) - \frac{1}{d+2}(1, 1, \dots,
1)$ in $\RR^{\hspace{.2mm}d+2}$. We know that the standard
$d$-sphere $S^{\,d}_{d+2}$ is the boundary complex of
$\Delta^{d+1}$. Let $\varphi : |S^{\,d}_{d + 2}| \to
S^{\,d+1}\subseteq \RR^{\hspace{.2mm}d+2}$ be the radial
projection. Then $S := \varphi(|S^{\,d}_{d + 2}|)$ is a
$d$-sphere in $\RR^{\hspace{.2mm}d+2}$ and $\varphi : |S^{\,d}_{d
+ 2}| \to S$ is a homeomorphism. Let $\alpha : S \to S$ be the
antipodal map. Then the quotient space $S/\alpha$ is the
$d$-dimensional real projective space $\RR P^{\,d}$.

Let ${\cal S}^{\hspace{.1mm}d} := (S^{\,d}_{d + 2})^{(1)}$ be the
first barycentric subdivision of $S^{\,d}_{d + 2}$. Let $V= \{v_1,
\dots, v_{d+2}\}$. Let ${\cal S}^{\hspace{.1mm}d}_a$ be the
abstract simplicial complex corresponding to ${\cal
S}^{\hspace{.1mm}d}$. We can identify $V({\cal
S}^{\hspace{.1mm}d}_a)$ with the set of proper subsets of $V$ by
$\widehat{\langle U \rangle} \mapsto U$. Let $\eta : V({\cal
S}^{\,d}_a) \to V({\cal S}^{\hspace{.1mm}d}_a)$ be given by
$\eta(U) = V \setminus U$. Then $\eta$ is an automorphism of
${\cal S}^{\hspace{.1mm}d}_a$ and $\eta \circ \eta = {\rm Id}$.
Observe that $d_{{\cal S}^{\hspace{.1mm}d}_a}(U, V\setminus U) =
3$ for any proper subset $U$ of $V$. (Since $d\geq 1$, $\#(U)$ or
$\#(V\setminus U)\geq 2$. Assume that $\#(U)\geq 2$ and $u\in U$.
Then $U\{u\}((V\setminus U)\cup \{u\}) (V\setminus U)$ is a path
of length 3 from $U$ to $\eta(U) = V\setminus U$.) So, $G = \{{\rm
Id}, \eta\}$ acts properly on ${\cal S}^{\hspace{.1mm}d}_a$.
Thus, by Proposition \ref{prop1}, ${\cal P}^{\,d} := {\cal
S}^{\hspace{.1mm}d}/G$ is a combinatorial $d$-manifold. Since the
number of vertices in ${\cal S}^{\hspace{.1mm}d}$ is $2^{d+2}
-2$, the number of vertices in ${\cal P}^{\hspace{.2mm}d}$ is
$2^{d+1} -1$.

Let $\eta_{gr} : V({\cal S}^{\hspace{.1mm}d}) \to V({\cal
S}^{\hspace{.1mm}d})$ be the simplicial map induced by $\eta$
(i.e., $\eta(\widehat{\langle U \rangle}) = \widehat{\langle V
\setminus U \rangle}$). Then $\eta_{gr}$ is an automorphism and
$\alpha\circ \varphi = \varphi \circ |\eta_{gr}|$. This implies
that $|{\cal S}^{\,d}/G| = |{\cal S}^{\,d}|/\eta_{gr}$ is
homeomorphic to $S/\alpha$. Thus, ${\cal P}^{\,d}$ is a
triangulation of $\RR P^{\,d}$. An explicit description of ${\cal
P}^{\,d}$ is given by Mukherjea in \cite{km}.}
\end{eg}

\begin{eg}$\!\!${\bf .} \label{eg21}
{\rm Consider the isometry group $G := \ZZ^d\!:\!S_d$ of
$\RR^{\hspace{.1mm}d}$, where the symmetric group $S_d$ acts on
$\RR^{\hspace{.1mm}d}$ by $(g(x))_i = x_{g(i)}$ for $1 \leq i \leq
d$ and $\ZZ^d$ acts by translations. Let $\sigma := \{(x_1,
\dots, x_d) : 0 \leq x_d \leq x_{d-1} \leq \cdots \leq x_1\leq
1\} \subseteq \RR^{\hspace{.1mm}d}$. Then $\sigma$ is a
$d$-simplex with vertices $(0, \dots, 0), (1, 0, \dots, 0), (1, 1,
0, \dots, 0), \dots, (1, \dots, 1, 0), (1, \dots, 1)$. Observe
that $\sigmaint \cap g(\sigma) = \emptyset$ for $1\neq g\in G$.
Consider the pure $d$-dimensional simplicial complexes
$$
Y = \{g(\sigma) : g \in S_d\} \mbox{ and } X = \{g(\sigma) : g
\in G\}.
$$
Then $Y$ triangulates the $d$-cube $[0, 1]^d$ and hence $X$
triangulates $\RR^{\hspace{.1mm}d}$. Clearly, $G$ is a group of
automorphisms of $X_a$ (the abstract simplicial complex
corresponding to $X$).

Let ${\cal L}_d := \{(x_1, \dots, x_{d}) : \sum_{i=1}^d
2^{i-1}x_i \equiv 0 \mbox{ mod } 2^{d+1} -1\}\subseteq \ZZ^d$.
Then ${\cal L}_d$ is a subgroup of $G$. (${\cal L}_d$ is a sub
lattice of $\ZZ^d$ and $\{(2^{d+1}-1, 0, \dots, 0), \dots, (-1,
0, \dots, 0, 2^{d+1-i}, 0, \dots, 0), \dots$, $(-1, 0, \dots, 0,
4)\}$ is a basis of ${\cal L}_d$.) Then ${\cal L}_d$ acts
properly on $X_a$. Therefore, by Proposition \ref{prop1}, $X/{\cal
L}_d$ is a combinatorial $d$-manifold. Since $[\ZZ^d : {\cal
L}_d] = 2^{d+1}-1$, it follows that the number of vertices in
$X/{\cal L}_d$ is $2^{d+1}-1$. The combinatorial $d$-manifold
$X/{\cal L}_d$ was constructed by K\"{u}hnel and La{\ss}mann in
\cite{kl3}. They have shown that $X/{\cal L}_d$ triangulates the
$d$-dimensional torus $S^{1}\times \cdots\times S^{1}$.}
\end{eg}

\section{Some General Results on Triangulations.}

In this section, we are presenting some results on triangulations.
Some of them are interesting and classical and some of them are
very useful.

\begin{theo} \label{DSE} $(${\bf Dehn-Sommerville Equations}$)${\bf
.} If $M$ is a combinatorial $d$-manifold then the $f$-vector and
the $h$-vector of $M$ satisfy the following.
\begin{enumerate}
\item[{\rm (i)}]
$\sum_{i= 0}^{d}(-1)^i f_i(M)  =  \chi(M)$ $(=0$ if $d$ is odd$)$.
\item[{\rm (ii)}]
If $d$ is even then $\sum_{i=2j-1}^{d}(-1)^i{i+1\choose
2j-1}f_i(M) = 0$  for $1\leq j\leq \frac{d}{2}$.
\item[{\rm (iii)}]
If $d$ is odd then $\sum_{i=2j}^{d}(-1)^i{i+1\choose 2j}f_i(M) =
0$ for $1\leq j\leq \frac{d-1}{2}$.
\item[{\rm (iv)}]
If $d = 2k$ then $h_j(M) - h_{d + 1 - j}(M) =(-1)^{d + 1 - j}{d +
1 \choose j}(\chi(M) - 2)$ for $0 \leq j \leq k$.
\item[{\rm (v)}] If $d = 2k-1$ then $h_j(M)-h_{d+1-j}(M)  = 0$
for $0 \leq j \leq k - 1$.
\end{enumerate}
\end{theo}

\noindent {\bf Proof.} The first equation is the Euler equation.

If $d$ is even then the link of a $(2j - 2)$-simplex is a
triangulation of the odd dimensional sphere $S^{\hspace{.2mm}d -
2j - 1}$ and hence the Euler characteristic of the link of a $(2j
- 2)$-simplex is $0$. If we take the sum of the Euler equations
of the links of all $(2j-2)$-simplices, then we get the second
equation.

Similarly, we get the third equation by taking the sum of the
Euler equations of the links of all the $(2j-1)$-simplices if $d$
is odd.

The last two equations follow from first three and the definition
of $h$-vector. \hfill $\Box$

\bigskip

Thus, for a combinatorial $d$-manifold $M$, $f_d(M), \dots,
f_{{(d+1)}/{2}}(M)$ can be express in terms of $f_0(M), \dots,
f_{{(d-1)}/{2}}(M)$ if $d$ is odd  and $f_d(M), \dots,
f_{{d}/{2}}(M)$ can be express in terms of $\chi(M), f_0(M),
\dots, f_{{d}/{2}-1}(M)$ if $d$ is even. Since, by the last
equation in Theorem \ref{DSE}, $(h_0(M), \dots, h_{\lfloor(d+1)/2
\rfloor}(M))$ determines the $h$-vector of $M$, it follows that
the $f$-vector of $M$ is determined by $(h_0(M), \dots,
h_{\lfloor(d+1)/2 \rfloor}(M))$. See \cite{kl} for more.

\medskip

For $n\geq d+2$ and $d\geq 1$, let $\varphi_d(n, d+1)  := dn - (d
+ 2)(d - 1)$ and $\varphi_k(n, d+1) := {d+1 \choose k}n - {d+2
\choose k+1}k$ for $1\leq k\leq d-1$. In the definition of
stacked sphere, we have seen that $f_k(S) = \varphi_k(n, d+1)$
for any $n$-vertex stacked $d$-sphere $S$ and $k\geq 1$. In
\cite{b1, b2}, Barnette proved the following\,:

\begin{theo}  $(${\bf Lower Bound Theorem}$)$. \label{LBT}
If $M$ is an $n$-vertex polytopal $d$-sphere then \vspace{-1mm}
\begin{enumerate}
\item[$(i)$] $f_k(M) \geq \varphi_k(n, d+1)$ for $1 \leq k \leq
d$ and
\item[$(ii)$] for $d \geq 3$, $f_d(M) = \varphi_d(n, d+1)$ if and
only if $M$ is a stacked sphere.
\end{enumerate}
\end{theo}

In \cite{b4}, Barnette proved the following generalization of
Theorem \ref{LBT} $(i)$.

\begin{theo} $\!\!${\bf .} \label{barnette4}
If $M$ is an $n$-vertex connected closed triangulated manifold of
dimension $d \geq 2$, then $f_k(M) \geq \varphi_k(n, d + 1)$ for
$1 \leq k \leq d$.
\end{theo}

Towards the classification of all the $n$-vertex triangulated
$d$-manifolds $M$ for which $f_k(M) = \varphi_k(n, d+1)$,
McMullen, Perles and Walkup observed the following independently
(see \cite{b2, kl, mc2}).

\begin{theo} $\!\!${\bf .} \label{MPW-reduction} Let $M$ be an
$n$-vertex $d$-dimensional $(d\geq 1)$ simplicial complex, such
that $f_1({\rm lk}_M(\sigma)) \geq \varphi_1(\deg_M(\sigma), d -
i)$ for any $i$-simplex $\sigma$ in $M$ $(0\leq i \leq d-2)$.
 \vspace{-1mm}
\begin{enumerate}
\item[$(i)$] Then $f_k(M) \geq \varphi_k(n, d+1)$ for $1 \leq k \leq
d$.
\item[$(ii)$] Moreover, if $f_k(M) = \varphi_k(n, d+1)$ for some $k
\geq 1$ then $f_1(M) = \varphi_1(n, d+1)$.
\end{enumerate}
\end{theo}

In \cite{ka}, Kalai showed that for $d \geq 3$, the edge graph of
any connected triangulated $d$-manifold is ``generically $(d +
1)$-rigid'' in the sense of rigidity of frameworks. The case
$k=1$ of Theorem \ref{barnette4} is an immediate consequence of
Kalai's rigidity theorem. Kalai also succeeds in using his
rigidity theorem to prove the following\,:

\begin{theo}$\!\!${\bf .} \label{kalai}
Let $M$ be an $n$-vertex triangulated closed manifold of
dimension $d\geq 3$. If $f_{1}(M) = \varphi_{1}(n, d+1)$, then
$M$ is a stacked $d$-sphere.
\end{theo}

In \cite{w}, Walkup proved Theorem \ref{barnette4} for $d = 3, 4$
and Theorem \ref{kalai} for $d=3$. For $d=2$ one observes the
following\,: If $M$ is an $n$-vertex connected combinatorial
2-manifold of Euler characteristic $\chi(M)$ then $f_1(M) =
3n-3\chi(M)$ and $f_2(M) = 2n - 2\chi(M)$. For every connected
combinatorial 2-manifold $M$, $\chi(M)\leq 2$ and $\chi(M) = 2$
if and only if $M$ is a (polytopal) 2-sphere. Thus, (i) $f_i(M)
\geq \varphi_i(n, 3)$ for $1\leq i \leq 2$ and (ii) $f_i(M) =
\varphi_i(n, 3)$ for $i =1$ or 2 if and only if $M$ is a
(polytopal) combinatorial 2-sphere. From Theorems
\ref{MPW-reduction} and \ref{kalai} one gets\,:

\begin{theo}$\!\!${\bf .} \label{kalai2}
Let $M$ be an $n$-vertex triangulated closed manifold of
dimension $d\geq 3$. If $f_{k}(M) = \varphi_{k}(n, d+1)$ for some
$k \geq 1$ then $M$ is a stacked $d$-sphere.
\end{theo}

In \cite{bd7}, we have presented a short and self-contained proof
of the following\,:

\begin{theo}$\!\!${\bf .} \label{LBT-NPM}
For $d \geq 2$, any $n$-vertex $d$-dimensional normal
pseudomanifold has $\geq n(d + 1) - {d + 2 \choose 2}$ edges. For
$d \geq 3$, equality holds only for stacked spheres.
\end{theo}

Let $C^{\, d}_{n}$ be the polytopal $d$-sphere as in Example
$\ref{eg15}$. Then $C^{\,d}_n$ is a $\lfloor \frac{d + 1}{2}
\rfloor$-neighbourly combinatorial $d$-manifold and hence
$h_j(C^{\,d}_n) = {n - d - 2 + j \choose j}$ for all $j = 0,
\dots, \lfloor \frac{d + 1}{2} \rfloor$. In \cite{mc1}, McMullen
proved the following\,:

\begin{theo}$\!\!${\bf .} \label{McMullen}
Let $X$ be a triangulation of the sphere $S^{\hspace{.4mm}d}$
with $n$ vertices. Then
\begin{enumerate}
\item[{\rm (i)}] If $h_j(X) \leq {n - d - 2 + j \choose j}$ for all $j = 0,
\dots, \lfloor \frac{d + 1}{2} \rfloor$ then $f_i(X) \leq
f_i(C^{\, d}_{n})$ for all $i = 0, \dots, d$.
\item[{\rm (ii)}] If $X$ is a polytopal $d$-sphere then $h_j(X) \leq {n
- d - 2 + j \choose j}$ for all $j = 0, \dots, \lfloor \frac{d +
1}{2} \rfloor$.
\end{enumerate}
\end{theo}

In \cite{st}, Stanley proved the following `Upper Bound
Conjecture' by showing that $h_j(X) \leq {n - d - 2 + j \choose
j}$ for all $j = 0, \dots, \lfloor \frac{d + 1}{2} \rfloor$
whenever $X$ triangulates $S^{\,d}$.

\begin{theo} $(${\bf Upper Bound Theorem for Spheres}$)$.
\label{UBTa} Let $X$ be an $n$-vertex simplicial complex. If $X$
triangulates $S^{\,d}$ then $f_i(X)\leq f_i(C^{\,d}_{n})$ for $1
\leq i \leq d$.
\end{theo}

For a combinatorial $d$-sphere, we get the following from Theorem
\ref{UBTa}\,:

\begin{cor}$\!\!${\bf .} \label{re4.6}
Let $M$ be an $n$-vertex $k$-neighbourly $d$-dimensional
pseudomanifold. If $M$ triangulates the $d$-sphere $S^{\,d}$ and
$n \geq d+3$ then $k\leq \lfloor\frac{d+1}{2}\rfloor$.
\end{cor}

\noindent {\bf Proof.} Since $C^{\,d}_n$ is not $(\lfloor \frac{d
+ 1}{2} \rfloor + 1)$-neighbourly for all $d$ with $n \geq d+3$,
the corollary follows from Theorem \ref{UBTa}. \hfill $\Box$

\bigskip

For a combinatorial $d$-sphere, Corollary \ref{re4.6} also
follows from Theorem \ref{DSE}.

\medskip

Let ${\cal T}$ and $\CC P^{\hspace{.2mm}2}_9$ be as in Examples
\ref{eg3} and \ref{eg10} respectively. Then their $f$-vectors are
as follows\,: $f({\cal T}) = (7, 21, 14)$, $f(\CC
P^{\hspace{.2mm}2}_9) = (9, 36, 84, 90, 36)$. Since the
$f$-vectors of any $S^{\hspace{.2mm}2}_7$ and $C^{\,4}_9$ are
$(7, 15, 10)$ and $(9, 36, 74, 75, 30)$ respectively, it follows
that the Upper bound theorem is not true for all manifolds. In
\cite{n}, Novik prove proved the following generalizations of
Theorem \ref{UBTa}.

\begin{theo} $(${\bf UBT for odd-dimensional
Homology Manifolds}$)$. \label{UBTb} Let $X$ be an $n$-vertex $(2k
- 1)$-dimensional homology manifold. Then $f_i(X) \leq
f_i(C^{\,2k - 1}_{n})$ for $1 \leq i \leq 2k - 1$.
\end{theo}

\begin{theo}$\!\!${\bf .} \label{UBTc}
For $d$ even, let $X$ be an $n$-vertex $d$-dimensional homology
manifold. If either
\begin{enumerate}
\item[{\rm (i)}] $d \equiv 0$ $($mod $4)$ and $\chi(X) \leq 2$, or
\item[{\rm (ii)}] $d \equiv 2$ $($mod $4)$, $\chi(X) \geq 2$ and
$H_{d/2}(|X|; \ZZ) = 0$
\end{enumerate}
then $f_i(X) \leq f_i(C^{\,d}_{n})$ for $1 \leq i \leq d$.
\end{theo}

In \cite{li}, Lickorish presented a proof of the following\,:

\begin{theo}$\!\!${\bf .} \label{re4.9}
Two simplicial complexes are combinatorially equivalent if and
only if they are stellar equivalent.
\end{theo}

Clearly, if two pseudomanifolds are bistellar equivalent then
they are combinatorially equivalent. In \cite{p}, Pachner proved
the following (see \cite{li} for a proof).

\begin{theo}$\!\!${\bf .} \label{re4.10}
Two combinatorial manifolds are combinatorially equivalent if and
only if they are bistellar equivalent.
\end{theo}

\section{Some Results on Minimal Triangulations.}

In this section, we are presenting some results without proofs.
Proofs are available in the references given.

Let $K$ be an $n$-vertex combinatorial 2-manifold. If $(n, f_1,
f_2)$ is the $f$-vector then $2f_1 = 3f_2$ and $f_1 \leq {n
\choose 2}$. Thus, $\chi(K)  = n - f_1 + f_2 = n- \frac{1}{3}f_1
\geq n - \frac{1}{3} {n \choose 2} = \frac{7n - n^2}{6}$. This
implies that $n \geq \frac{1}{2}(7 + \sqrt{49 -24\chi(K)})$. It
is known that the Klein bottle (whose Euler characteristic is 0)
has an 8-vertex triangulation and has no 7-vertex triangulation
(Theorem \ref{d=2n=7} below). From the classification of 8-vertex
combinatorial 2-manifolds (Theorem \ref{d=2n=8} below), we know
that there is no 8-vertex combinatorial 2-manifold of Euler
characteristic $-1$. In two articles (\cite{r, jr}), Ringel and
Jungerman proved the following\,:

\begin{theo}$\!\!${\bf .} \label{th1}
Let $M$ be a closed surface which is not the Klein bottle, the
double torus or the non-orientable surface of Euler
characteristic $-1$. Then $M$ has an $n$-vertex triangulation if
and only if $n \geq \frac{1}{2}(7 + \sqrt{49 -24\chi(M)})$. In
each of those three cases, one needs one more vertex for
triangulations.
\end{theo}

It is known that the only combinatorial 2-manifolds on at most 6
vertices are $S^{\,2}_4$, $S^{\,0}_2 \ast S^{1}_3$, $S^{\,0}_2
\ast S^{\,0}_2 \ast S^{\,0}_2$, $\Sigma^1(S^{1}_5)$ and $\RR
P^{\,2}_6$. In \cite{d2}, we have shown the following\,:

\begin{theo}$\!\!${\bf .} \label{d=2n=7}
There are exactly nine $7$-vertex combinatorial $2$-manifolds,
five of which triangulate the $2$-sphere $S^{\hspace{.4mm}2}$,
three of which triangulate $\RR P^{\hspace{.2mm}2}$ and one
triangulate $S^{\hspace{.2mm}1}\times S^{\hspace{.2mm}1}$.
\end{theo}

In \cite{dn}, we have proved the following\,:

\begin{theo}$\!\!${\bf .} \label{d=2n=8}
There are exactly $44$ distinct combinatorial $2$-manifolds on
$8$ vertices. One of these combinatorial $2$-manifolds consists of
two copies of $S^{\,2}_4$'s,  $14$ of these triangulate
$S^{\,2}$, $16$ triangulate $\RR P^{\,2}$, seven triangulate
$S^{1}\times S^{1}$ and six triangulate the Klein bottle.
\end{theo}

For $g\geq 0$, let $M(g, +)$ denote the orientable surface of
genus $g$ and let $M(g, -)$ denote the non-orientable surface of
genus $g$. (So, $M(1, +) = S^{\hspace{.2mm}1} \times
S^{\hspace{.2mm}1}$ and $M(2, -)$ is the Klein bottle.) Thus,
$\chi(M(g, +)) = 2-2g$ and $\chi(M(g, -)) = 2 - g$. In \cite{lu5,
sl}, Lutz and Sulanke have enumerated (via computer search) all
the triangulated 2-manifolds with at most 12 vertices. They have
shown the following\,:

\begin{theo}$\!\!${\bf .} \label{d=2n=9}
There are precisely $655$ combinatorial $2$-manifold with $9$
vertices: $50$ of these triangulate $S^{\hspace{.2mm}2}$, $112$
triangulate $S^{\hspace{.1mm}1}\times S^{\hspace{.1mm}1}$, $134$
triangulate $\RR P^{\hspace{.2mm}2}$, $187$ triangulate the Klein
bottle, $133$ triangulate $M(3, -)$, $37$ triangulate $M(4, -)$
and $2$ triangulate $M(5, -)$.
\end{theo}

\begin{theo}$\!\!${\bf .} \label{d=2n=10}
There are precisely $42\hspace{.3mm}426$ combinatorial
$2$-manifold with $10$ vertices: $233$ of these triangulate
$S^{\hspace{.2mm}2}$, $2\hspace{.3mm}109$ triangulate
$S^{\hspace{.1mm}1}\times S^{\hspace{.1mm}1}$, $865$ triangulate
$M(2, +)$, $20$ triangulate $M(3, +)$, $1\hspace{.2mm}210$
triangulate $\RR P^{\hspace{.2mm}2}$, $4\hspace{.3mm}462$
triangulate the Klein bottle, $11784$ triangulate $M(3, -)$,
$13\hspace{.3mm}657$ triangulate $M(4, -)$, $7\hspace{.3mm}050$
triangulate $M(5, -)$, $1\hspace{.2mm}022$ triangulate $M(6, -)$
and $14$ triangulate $M(7, -)$.
\end{theo}

\begin{theo}$\!\!${\bf .} \label{d=2n=11}
There are precisely $11\hspace{.2mm}590\hspace{.3mm}894$
combinatorial $2$-manifold with $11$ vertices: $1\hspace{.2mm}249$
of these triangulate $S^{\hspace{.2mm}2}$, $37\hspace{.2mm}867$
triangulate $S^{\hspace{.1mm}1} \times S^{\hspace{.1mm}1}$,
$113\hspace{.2mm}506$ triangulate $M(2, +)$, $65\hspace{.3mm}876$
triangulate $M(3, +)$, $821$ triangulate $M(4, +)$,
$11\hspace{.2mm}719$ triangulate $\RR P^{\hspace{.2mm}2}$,
$86\hspace{.3mm}968$ triangulate the Klein bottle,
$530\hspace{.3mm}278$ triangulate $M(3, -)$,
$1\hspace{.2mm}628\hspace{.3mm}504$ triangulate $M(4, -)$,
$3\hspace{.3mm}355\hspace{.3mm}250$ triangulate $M(5, -)$,
$3\hspace{.3mm}623\hspace{.3mm}421$ triangulate $M(6, -)$,
$1\hspace{.2mm}834\hspace{.2mm}160$ triangulate $M(7, -)$,
$295\hspace{.3mm}291$ triangulate $M(7, -)$ and
$5\hspace{.3mm}982$ triangulate $M(9, -)$.
\end{theo}

\begin{theo}$\!\!${\bf .} \label{d=2n=12}
There are precisely $12\hspace{.3mm}561\hspace{.3mm}206
\hspace{.3mm}794$ combinatorial $2$-manifold with $12$ vertices:
$7,595$ of these triangulate $S^{\hspace{.2mm}2}$,
$605\hspace{.3mm}496$ triangulate $S^{\hspace{.1mm}1}\times
S^{\hspace{.1mm}1}$, $7\hspace{.3mm}085\hspace{.3mm}444$
triangulate $M(2, +)$, $25\hspace{.3mm}608\hspace{.3mm}643$
triangulate $M(3, +)$, $14\hspace{.3mm}846\hspace{.3mm}522$
triangulate $M(4, +)$, $751\hspace{.2mm}593$ triangulate $M(5,
+)$, $59$ triangulate $M(6, +)$, $114\hspace{.3mm}478$
triangulate $\RR P^{\hspace{.2mm}2}$, $1\hspace{.2mm}448
\hspace{.3mm}516$ triangulate the Klein bottle, $16\hspace{.3mm}
306\hspace{.3mm}649$ triangulate $M(3, -)$, $99\hspace{.3mm}
694\hspace{.3mm}693$ triangulate $M(4, -)$, $473\hspace{.3mm}
864\hspace{.3mm}807$ triangulate $M(5, -)$, $1\hspace{.2mm}
479\hspace{.2mm}135\hspace{.3mm}833$ triangulate $M(6, -)$,
$3\hspace{.2mm}117\hspace{.3mm}091\hspace{.2mm}975$ triangulate
$M(7, -)$, $3\hspace{.3mm}935\hspace{.3mm}668\hspace{.3mm}832$
triangulate $M(8, -)$, $2\hspace{.3mm}627\hspace{.3mm}619
\hspace{.3mm}810$ triangulate $M(9, -)$, $711\hspace{.2mm}
868\hspace{.3mm}010$ triangulate $M(10, -)$, $49\hspace{.3mm}
305\hspace{.3mm}639$ triangulate $M(11, -)$ and $182\hspace{.3mm}
200$ triangulate $M(12, -)$.
\end{theo}

We know (see Theorem \ref{n=d+4b} below) that a combinatorial
3-manifold on at most 7 vertices is a polytopal 3-sphere. In
\cite{a1}, Altshuler proved the following\,:

\begin{theo}$\!\!${\bf .} \label{d=3n=8a}
Every combinatorial $3$-manifold with at most $8$ vertices is a
combinatorial $3$-sphere.
\end{theo}

In \cite{gs}, Gr\"{u}nbaum and Sreedharan shown the following\,:

\begin{theo}$\!\!${\bf .} \label{d=3n=8b}
There are exactly $37$ polytopal $3$-spheres on $8$ vertices.
\end{theo}

Gr\"{u}nbaum and Sreedharan have also constructed the $8$-vertex
non-polytopal $3$-sphere $S^{\, 3}_{8, 38}$ (see Example
\ref{eg5}). In \cite{b3}, Barnette have proved the following\,:

\begin{theo}$\!\!${\bf .} \label{d=3n=8c}
There are exactly two non-polytopal combinatorial $3$-sphere on
$8$ vertices, namely, $S^{\, 3}_{8, 38}$ and $S^{\, 3}_{8, 39}$
$($given in Example $\ref{eg5})$.
\end{theo}

So, there are exactly 39 combinatorial 3-manifolds with 8
vertices. We got a different proof of this. This follows from the
next two theorems which we have proved in \cite{dn4}.

\begin{theo}$\!\!${\bf .} \label{d=3n=8d}
Every $8$-vertex $3$-pseudomanifold is obtained from a neighbourly
$8$-vertex $3$-pseudomanifold by a sequence of bistellar
$2$-moves.
\end{theo}

\begin{theo}$\!\!${\bf .} \label{d=3n=8e}
If $M$ is an $8$-vertex neighbourly  combinatorial $3$-manifold
then $M$ is isomorphic to one of $S^{\,3}_{8, 35}, \dots,
S^{\,3}_{8, 38}$ $($given in Example $\ref{eg5})$.
\end{theo}

In Example \ref{eg16}, we have seen that there exists a 9-vertex
triangulation (namely, $K^{\hspace{.2mm}3}_9$) of the twisted
product $\TPSS$ and there exists a 10-vertex triangulation
(namely, $K^{\hspace{.2mm}3}_{10}$) of $S^{\hspace{.2mm}2} \times
S^{\hspace{.2mm}1}$. In \cite{w}, Walkup proved the following\,:

\begin{theo}$\!\!${\bf .} \label{walkup1}
There exists an $n$-vertex triangulation of $~ S^{\,2} \times
S^{\hspace{.2mm}1}$ only if $n \geq 10$.
\end{theo}

\begin{theo}$\!\!${\bf .} \label{walkup2}
If $K$ is a combinatorial $3$-manifold and $|K|$ is not
homeomorphic to $S^{\,3}$, $\TPSS$ or $S^{\,2}\times S^1$ then
$f_1(K) \geq 4f_0(K) +8$ and hence $f_0(K) \geq 11$.
\end{theo}

Thus, for a combinatorial triangulation of $\RR
P^{\hspace{.2mm}3}$ one needs at least 11 vertices. Therefore,
from Example \ref{eg6} and Theorem \ref{walkup2} one gets\,:

\begin{cor}$\!\!${\bf .} \label{walkup3}
There exists an $n$-vertex triangulation of $~ \RR
P^{\hspace{.2mm}3}$ if and only if $n\geq 11$.
\end{cor}

In \cite{as1, as2}, Altshuler and Steinberg showed (via a computer
search) the following\,:

\begin{theo}$\!\!${\bf .} \label{d=3n=9a}
There are exactly $1297$ combinatorial $3$-manifolds on nine
vertices. One of these is $K^3_9$ and other $1296$ are
combinatorial $3$-spheres. Among these $1296$ combinatorial
spheres, $50$ are neighbourly. Among these $50$ neighbourly
combinatorial spheres, $23$ are polytopal and $27$ are
non-polytopal.
\end{theo}

Altshuler and Steinberg also showed (using computer) that any two
of these 1296 spheres are bistellar equivalent via a finite
sequence of proper bistellar moves. In \cite{bd6}, we present
computer-free proofs of the following\,:

\begin{theo}$\!\!${\bf .} \label{d=3n=9b}
Every $9$-vertex combinatorial $3$-manifold is obtained from a
neighbourly $9$-vertex combinatorial $3$-manifold by a sequence
of $($at most $10)$ bistellar $2$-moves.
\end{theo}

\begin{theo}$\!\!${\bf .} \label{d=3n=9c}
Up to isomorphism, there is a unique $9$-vertex non-sphere
combinatorial $3$-manifold, namely $K^{\hspace{.2mm}3}_9$.
\end{theo}

In \cite{lu4, sl}, Lutz and Sulanke have enumerated (via computer
search) all the triangulated 3-manifolds with 10 and 11 vertices.
They have shown the following\,:

\begin{theo}$\!\!${\bf .} \label{d=3n=10}
There are precisely $249\hspace{.3mm}015$ combinatorial
$3$-manifold with $10$ vertices: $247\hspace{.3mm}882$ of these
triangulate the $3$-sphere $S^{\hspace{.2mm}3}$, $615$
triangulate the twisted product $\TPSS$ and $518$ triangulate the
sphere product $S^{\hspace{.2mm}2} \times S^{\hspace{.2mm}1}$.
\end{theo}

\begin{theo}$\!\!${\bf .} \label{d=3n=11}
There are precisely $172\hspace{.3mm}638\hspace{.3mm}650$
combinatorial $3$-manifolds with $11$ vertices: $166\hspace{.3mm}
564\hspace{.3mm}303$ of these triangulate the $3$-sphere
$S^{\hspace{.2mm}3}$, $3\hspace{.2mm}116\hspace{.3mm}818$
triangulate the twisted sphere product $\TPSS$, $2\hspace{.3mm}
957\hspace{.3mm}499$ triangulate the sphere product
$S^{\hspace{.2mm}2} \times S^{\hspace{.2mm}1}$ and $30$
triangulate the real projective $3$-space $\RR
P^{\hspace{.2mm}3}$.
\end{theo}

To get an estimate of the minimal number of vertices for a
triangulation of a 4-manifold in terms of the Euler
characteristic, K\"{u}hnel has proved the following (in
\cite{k1})\,:

\begin{theo}$\!\!${\bf .} \label{d=4a}
If $M$ is a combinatorial $4$-manifold with $n$ vertices then
$10(\chi(M) - 2) \leq {n - 4 \choose 3}$. Equality holds if and
only if $M$ is $3$-neighbourly.
\end{theo}

Since the Euler characteristic of the K3 surface is 24, by Theorem
\ref{d=4a}, any combinatorial triangulation of the K3 surface
requires 16 vertices. In \cite{ck}, Casella and K\"{u}hnel have
constructed a 16-vertex triangulation of the K3 surface (${\rm
K3}_{16}$ in Example \ref{eg13}). It follows from Theorem
\ref{d=4a} that any combinatorial triangulation of
$(S^{\,2}\times S^{\,2}) \# (S^{\,2}\times S^{\,2})$ requires 12
vertices. In \cite{lu1}, Lutz has constructed two
(non-isomorphic) 12-vertex combinatorial triangulations of
$(S^{\,2}\times S^{\,2}) \# (S^{\,2}\times S^{\,2})$.

Observe that the equality holds in Theorem \ref{d=4a} for
$S^{\hspace{.3mm}4}_6$, $\CC P^{\hspace{.2mm}2}_9$ and ${\rm
K3}_{16}$. In \cite{kl1}, K\"{u}hnel and La{\ss}mann proved the
following\,:

\begin{theo}$\!\!${\bf .} \label{d=4b}
Let $M$ be an $n$-vertex combinatorial $4$-manifold. If $n \leq
13$ and $M$ is $3$-neighbourly then $M = S^{\hspace{.3mm}4}_6$ or
$\CC P^{\hspace{.2mm}2}_9$.
\end{theo}

For negative Euler characteristic, we get a lower bound of number
of vertices from the following result of Walkup \cite{w}\,:

\begin{theo}$\!\!${\bf .} \label{d=4c}
If $M$ is an $n$-vertex combinatorial $4$-manifold then $f_1(M)
\geq 5n - \frac{15}{2}\chi(M)$. Equality holds if and only if the
links of all the vertices are stacked $3$-spheres.
\end{theo}

Since $f_1(M) \leq {n \choose 2}$ for any $n$-vertex simplicial
complex $M$, from Theorem \ref{d=4c}, we get the following\,:

\begin{cor}$\!\!${\bf .} \label{d=4d}
If $M$ is an $n$-vertex combinatorial $4$-manifold then $-15
\chi(M) \leq n(n-11)$. Equality implies $M$ is $2$-neighbourly.
\end{cor}

A $d$-dimensional pseudomanifold has at least $d+2$ vertices. It
is easy to see that the only $d$-pseudomanifold with $d+2$
vertices is $S^{\hspace{.3mm}d}_{d+2}$. It is also known that a
combinatorial $d$-sphere on $d+3$ vertices is a join of standard
spheres. In \cite{bd2}, we have seen the following\,:

\begin{theo}$\!\!${\bf .} \label{n=d+3}
If $M$ is a $d$-dimensional $(d \geq 1)$ pseudomanifold with $d +
3$ vertices then $M$ is a polytopal sphere and is isomorphic to
$S^{\,c}_{c + 2} \ast S^{\hspace{.3mm}d - c - 1}_{d - c + 1}$ for
some $c < d$.
\end{theo}

Thus, $(d + 3)$-vertex $d$-dimensional pseudomanifolds are
completely reducible. In \cite{m}, Mani has proved the
following\,:

\begin{theo}$\!\!${\bf .} \label{n=d+4a}
Every combinatorial $d$-spheres on at most $d + 4$ vertices is
polytopal.
\end{theo}

In \cite{bd2}, we have classified all the $d$-dimensional
pseudomanifold on $d + 4$ vertices.  In particular, we have proved
the following\,:

\begin{theo}$\!\!${\bf .} \label{n=d+4b}
For $n \geq 6$, the $n$-vertex combinatorial $(n-4)$-manifolds
consist of\,:
\begin{enumerate}
\item[$({\rm a})$] The $6$-vertex combinatorial $2$-manifold $\RR
P^2_6$ $($defined in Example $\ref{eg2})$,
\item[$({\rm b})$] completely reducible polytopal spheres; their
number is $ \left\lfloor\frac{n(n-6)}{12}\right\rfloor +1$, and
\item[$({\rm c})$] irreducible polytopal spheres; their number is
the integer nearest to
$$
2^{ \left\lfloor (n-3)/2 \right\rfloor} - \frac{1}{12} n^{2} - 1 +
\frac{1}{4n} \sum_{r} \varphi (r) 2^{n/r}\ .
$$
Here $\varphi$ is Euler's totient function and the sum is over all
the odd divisors $r$ of $n$.
\end{enumerate}
\end{theo}

In \cite{bg}, Barnette and Gannon proved the following\,:

\begin{theo}$\!\!${\bf .} \label{n<d+6}
Let $M$ be an $n$-vertex combinatorial $d$-manifold, where $d \geq
3$ and $n \leq d + 5$. If $d\neq 4$ then $M$ is a combinatorial
$d$-sphere.
\end{theo}

In \cite{bk}, Brehm and K\"{u}hnel proved the following more
general results\,:

\begin{theo}$\!\!${\bf .} \label{BrehmKuhnel}
Let $M^{\,d}_n$ be an $n$-vertex combinatorial $d$-manifold
$(d>0)$. \vspace{-2mm}
\begin{description}
\item[(a)]  If $n < 3\lceil d/2 \rceil + 3$ then $M^{\,d}_n \approx
S^{\,d}_{d + 2}$.
\item[(b)] If $n = 3d/2 + 3$ and $M^{\,d}_n \not\approx S^{\,d}_{d
+ 2}$ then $d = 2, 4, 8$ or $16$. Moreover, $M^{\hspace{.2mm}2}_6
= \RR P^{\hspace{.2mm}2}_6$, $M^{\,4}_9$ triangulates $\CC
P^{\hspace{.2mm}2}$ and for $d = 8$ or $16$,
$|M^{\hspace{.2mm}d}_n|$ is a simply connected cohomology
projective plane over quaternions or Cayley numbers, respectively.
\end{description}
\end{theo}

In \cite{kb}, K\"{u}hnel and Banchoff constructed a 9-vertex
triangulation of $\CC P^{\hspace{.2mm}2}$ (see Example
\ref{eg10}). In \cite{kl1}, K\"{u}hnel and La{\ss}mann showed (by
the help of a computer) the following\,:

\begin{theo}$\!\!${\bf .} \label{cp29}
Up to isomorphism there is a unique $9$-vertex triangulation of
$~\CC P^{\hspace{.2mm}2}$.
\end{theo}

Computer-free proofs of the uniqueness of $\CC P^2_9$ have
appeared in \cite{am} and \cite{bd1}. In \cite{bd3}, we have
presented a very short (theoretical) proof of the uniqueness of
$\CC P^{\,2}_9$.

In \cite{bk2}, Brehm and K\"{u}hnel constructed three 15-vertex
combinatorial 8-manifolds of Euler characteristic 3. They also
showed that these three triangulate the same pl manifold. So, we
have\,:

\begin{theo}$\!\!${\bf .} \label{hp215}
There exist at least three different $15$-vertex combinatorial
$8$-manifolds which are not combinatorial spheres.
\end{theo}

In \cite{am}, Arnoux and Marin proved the following\,:

\begin{theo}$\!\!${\bf .} \label{n=3d/2+3a}
If $M$ is a non-sphere combinatorial $d$-manifold on $3d/2 + 3$
vertices then $M$ satisfies complementarity.
\end{theo}

In \cite{d1}, we have proved the following converse\,:

\begin{theo}$\!\!${\bf .} \label{n=3d/2+3b}
Let $M$ be an $n$-vertex combinatorial $d$-manifold. If $M$
satisfies complementarity then $d = 2, 4, 8$ or $16$ with $n =
3d/2 + 3$ and $M$ is a non-sphere.
\end{theo}

In \cite{d3, bd4}, we have shown the following\,:

\begin{theo}$\!\!${\bf .} \label{comp}
Let $M$ be an $n$-vertex $d$-dimensional pseudomanifold with
complementarity. If $n \leq d+6$ or $d\leq 6$ then $M$ is either
$\RR P^{\hspace{.2mm}2}_6$ or $\CC P^{\hspace{.2mm}2}_9$.
\end{theo}

We know from Theorem \ref{walkup2} that the minimal number of
vertices required for a triangulation of $\RR P^{\hspace{.2mm}3}$
is 11. We have seen in Example \ref{eg12} that there exists a
16-vertex triangulation of $\RR P^{\hspace{.3mm}4}$. In \cite{am},
Arnoux and Marin proved the following\,:

\begin{theo}$\!\!${\bf .} \label{rpdn}
Let $M$ be an $n$-vertex combinatorial $d$-manifold. If the
cohomology ring of $|M|$ is same as that of $~\RR
P^{\hspace{.3mm}d}$ then $n \geq {d + 2 \choose 2}$. Moreover,
equality is possible only for $d = 1$ and $d = 2$.
\end{theo}

\begin{theo}$\!\!${\bf .} \label{cpdn}
Let $M$ be an $n$-vertex combinatorial $2d$-manifold. If the
cohomology ring of $|M|$ is same as that of $~\CC
P^{\hspace{.3mm}d}$ then $n \geq (d+1)^2$. Moreover, equality is
possible only for $d = 1$ and $d = 2$.
\end{theo}

In \cite{bk}, Brehm and K\"{u}hnel proved the following\,:

\begin{theo}$\!\!${\bf .} \label{n<2d+4+i}
Let $M$ be a combinatorial $d$-manifold with $n$ vertices. If $n
\leq 2d + 3-i$ for some $i$ with $1\leq i < d/2$ then $|M|$ is
$i$-connected.
\end{theo}

Thus, if $m \geq n \geq 1$ then for a combinatorial triangulation
of $S^{\hspace{.4mm} m} \times S^{\hspace{.4mm} n}$ we need at
least $2m + n + 4$ vertices. In Example \ref{eg14}, we have seen
that there exists a combinatorial triangulation of $S^{\,3} \times
S^{\,2}$ with 12 vertices. In \cite{lu1}, Lutz has proved the
following\,:

\begin{theo}$\!\!${\bf .} \label{s3xs3}
There are at least two $13$-vertex combinatorial triangulations
of $S^{\,3}\times S^{\,3}$.
\end{theo}

In \cite{bd5}, we have proved the following\,:

\begin{theo}$\!\!${\bf .} \label{z2hsa}
Let $M$ be an $n$-vertex combinatorial $d$-manifold. If $|M|$ is a
$\ZZ_2$-homology sphere and $n \leq d+8$ then $M$ is a
combinatorial sphere.
\end{theo}

\begin{theo}$\!\!${\bf .} \label{z2hsb}
Let $M$ be a $(d + 9)$-vertex combinatorial triangulation of a
$\ZZ_2$-homology $d$-sphere. If $M$ is not a combinatorial sphere
then $M$ can not admit any bistellar $i$-move for $i < d$.
\end{theo}

We have seen in Example \ref{eg7} that there exists a 12-vertex
combinatorial triangulation of the lens space $L(3, 1)$. Since
$L(3, 1)$ is a $\ZZ_2$-homology $3$-sphere, Theorem \ref{z2hsa} is
sharp for $d = 3$. It follows from Theorem \ref{z2hsb} that a
12-vertex combinatorial triangulation of $L(3, 1)$ can not admit
any bistellar $i$-move for $0 \leq i \leq 2$.

From Theorem \ref{n<2d+4+i}, we know that a triangulation of a
closed pl manifold of dimension $d \geq 3$ requires at least $2d
+ 3$ vertices. We also know that there exist such triangulations
(namely, K\"{u}hnel's combinatorial $d$-manifold $K^{d}_{2d + 3}$
in Example \ref{eg16}) with $(2d + 3)$ vertices. In \cite{bd7} we
have proved the following\,:

\begin{theo}$\!\!${\bf .} \label{n=2d+3}
For $d \geq 3$, K\"{u}hnel's complex $K^{d}_{2d+3}$ is the only
non-simply connected $(2d+3)$-vertex triangulated manifold of
dimension $d$.
\end{theo}

\noindent {\bf Acknowledgement\,:} The author would like to thanks
Ashish Kumar Upadhyay, who read an earlier version of the
manuscript and made important suggestions. The author wishes to
thank Bhaskar Bagchi for making him aware of several important
papers in the literature. The author gratefully acknowledges the
financial support obtained from DST (Grant: SR/S4/MS-272/05) and
from UGC-SAP/DSA-IV.



\end{document}